\documentclass[natbib]{svjour3}

\smartqed  

\usepackage{amsmath,amssymb,graphicx,xcolor,url,caption,float}

\RequirePackage{fix-cm}

\newcommand{\uvec}{\mathbf{u}}

\newcommand{\wvec}{\mathbf{w}}
\newcommand{\xvec}{\mathbf{x}}
\newcommand{\yvec}{\mathbf{y}}
\newcommand{\zvec}{\mathbf{z}}

\newcommand{\Svec}{\mathbf{S}}
\newcommand{\Tvec}{\mathbf{T}}
\newcommand{\Rvec}{\mathbf{R}}
\newcommand{\Xvec}{\mathbf{X}}
\newcommand{\Yvec}{\mathbf{Y}}
\newcommand{\Zvec}{\mathbf{Z}}

\newcommand{\pr}{\mathrm{Pr}}
\newcommand{\E}{\mathrm{E}}
\newcommand{\var}{\mathrm{Var}}
\newcommand{\cov}{\mathrm{Cov}}
\newcommand{\cor}{\mathrm{Cor}}
\newcommand{\Xd}{X_1,X_2,\cdots,X_d}
\newcommand{\ks}{k_\sigma}
\newcommand{\cb}{I_{\sigma}(\mathbf{X})}
\newcommand{\scb}{I^2_{\sigma}(\mathbf{X})}
\newcommand{\csd}{C_{\sigma,d}}
\newcommand{\ecb}{\widehat{I}_{\sigma,n}(\mathbf{X})}
\newcommand{\secb}{\widehat{I}_{\sigma,n}^2(\mathbf{X})}
\newcommand{\cx}{\mathsf{C}_{\Xvec}}
\newcommand{\cxn}{\mathsf{C}_{\xvec,n}}
\newcommand{\kss}{k_{\frac{\sigma}{\sqrt{2}}}}
\newcommand{\Rd}{\mathbb{R}^d}

\newcommand{\intc}{\int_{[0,1]^d}}

\newtheorem{innercustomgeneric}{\customgenericname}
\providecommand{\customgenericname}{}
\newcommand{\newcustomtheorem}[2]{%
	\newenvironment{#1}[1]
	{%
		\renewcommand\customgenericname{#2}%
		\renewcommand\theinnercustomgeneric{##1}%
		\innercustomgeneric
	}
	{\endinnercustomgeneric}
}

\newcustomtheorem{customthm}{Theorem}
\newcustomtheorem{customlemma}{Lemma}

\begin{document}

\title{Some New Copula Based Distribution-free Tests of Independence among Several Random Variables}

\titlerunning{Copula Based Distribution-Free Tests}

\author{Angshuman Roy \and Anil K. Ghosh \and Alok Goswami \and C. A. Murthy}

\institute{A. Roy \at
              Indian Statistical Institute\\
              Applied Statistics Unit\\
              \email{angshuman.roy.89@gmail.com}
           \and
           A. K. Ghosh \at
              Indian Statistical Institute\\
              Stat Math Unit\\
              \email{akghosh@isical.ac.in}
           \and
           A. Goswami \at
	           Indian Statistical Institute\\
	           Stat Math Unit\\
	           \email{alok@isical.ac.in}
	       \and
	       C. A. Murthy \at
	       Indian Statistical Institute\\
	       Machine Intelligence Unit\\
	       \email{murthy@isical.ac.in}
}
\date{Received: date / Accepted: date}

\maketitle

\begin{abstract}
Over the last couple of decades, several copula based methods have been proposed in the literature to test for independence among several random variables. But these existing tests are not invariant under monotone transformations of the variables, and they often perform poorly if the dependence among the variables is highly non-monotone in nature. In this article, we propose a copula based measure of dependency and use it to construct some distribution-free tests of independence. The proposed measure and the resulting tests, all are invariant under permutations and strictly monotone transformations of the variables. Our dependency measure involves a kernel function with an associated bandwidth parameter. 
We adopt a multi-scale approach, where we look at the results obtained for several choices of the bandwidth 
and aggregate them judiciously. Large sample properties of the dependency measure and the resulting tests are derived under appropriate regularity conditions. Several simulated and real data sets are analyzed to compare the performance of the proposed tests with some popular tests available in the literature.
\end{abstract}

\section{Introduction}

Measuring and testing for dependence among $d~(d\ge 2)$ random variables is a classical problem in statistics. For $d=2$, Pearson's product-moment correlation coefficient, being arguably the most popular measure of dependence, has been used to construct test of independence between two random variables \citep[see, e.g.,][]{anderson2003}. But this measure is sensitive against outliers and extreme values, and it often fails to capture non-linear dependence between the variables. Other popular measures of association like Spearman's rank correlation coefficient $\rho$ \citep{spearman1904}, Kendall's measure of association $\tau$ \citep{kendall1938}, Blomqvist's quadrant statistic $\beta$ \citep{blomqvist1950}, are robust against outliers. They can also detect monotone or near monotone relationship between the variables. Tests based on these rank based statistics have the distribution-free property under the null hypothesis of independence. \cite{hoeffding1948} also developed a distribution-free test of independence using  $\phi$ statistic based on bivariate copula. \cite{reshef2011detecting} proposed maximal information coefficient as a measure of dependency, but the tests based on this measure usually have low powers. Tests of independence between two random vectors include the work of  \cite{gieser1997,taskinen2003,taskinen2005, heller2012,heller2013,biswas2016,sarkar2018some}. \cite{szekely2007measuring} developed a test of independence based on distance correlation, which is known as the dCov test .  \cite{gretton2007kernel} considered a test based on the Hilbert-Schmidt norm of the cross-covariance operator, which is popularly known as the HSIC test.

\cite{um2001} generalized Giesser and Randles' (1997) test  for several variables. \cite{gaisser2010multivariate} developed multivariate extensions of Hoeffding's (1948) $\phi$ statistic and the associated test. \cite{pfister2016kernel} and \cite{fan2017} proposed multivariate extensions of the HSIC test and the dCov test (referred to as the dHSIC test and the mdCov test), respectively. 

Using the ideas of copula, \cite{nelsen1996nonparametric} and \cite{flores05} proposed multivariate extensions of Spearman's $\rho$, Kendall's $\tau$ and Blomqvist's $\beta$ statistics. Tests of independence based on these dependency measures have the distribution-free property, but they often yield poor performance when the relationships among the variables are highly non-monotone in nature. To take care of this problem, in this article, we propose a new copula based multivariate measure of dependency and use it to test mutual independence among several random variables.

Our work is motivated by \cite{poczos2012}, where the authors proposed a dependency measure, which is non-negative and takes the value $0$ if and only if the random variables $X_1,X_2,\ldots,X_d$ are jointly independent. This measure is also invariant under permutation of the coordinates. Our proposed measure satisfies all these desirable properties. Moreover, while the dependency measure of \cite{poczos2012} is only invariant under strictly increasing transformations of the $X_i$'s, our measure is invariant under strictly monotone transformations. 
It also satisfies several other desirable properties. 
For instance, in the bivariate case, it satisfies all Dependence Axioms proposed by \cite{schweizer1981nonparametric}. We propose a data driven estimator 
of this measure, which also enjoys similar desirable properties. 
Unlike \cite{poczos2012}, here one does not need to generate 
observations from a uniform distribution for its construction. 
Our dependency measure and its estimator involve a kernel function, and we use the 
Gaussian kernel for this purpose.

One can use this estimator to a develop distribution-free test. 
However, the performance of the test depends on the bandwidth parameter associated with the Gaussian kernel. While larger bandwidths work well for near monotone relationships (i.e., when the conditional expectation of one variable is nearly a monotone function of others) among the variables, smaller bandwidths are preferred to detect non-monotone relations. So, borrowing idea from multi-scale  classification \citep[see, e.g.,][]{ghosh2006}, here we adopt a multi-scale approach, where we look at the results for various choices of bandwidth and then aggregate them judiciously to arrive at the final decision.

The organization of the paper is as follows. In Section 2, we define our copula based dependency measure and derive some of its desirable properties. In particular, we prove its invariance over permutations and monotone transformations of the variables. In Section 3, we propose a nonparametric estimate of this dependency measure and investigate its theoretical properties. Some distribution-free tests based on this estimate are constructed in Section 4, where we also establish the large sample consistency of these tests. Several simulated and real data sets are analyzed in Section 5 to compare the performance of the proposed tests with some popular tests available in the literature. Section 6 contains a brief summary of the work and ends with a discussion on some possible directions for future research. All proofs and mathematical details are given in Appendix Section.

\section{The proposed measure of dependency}

Our measure of dependency is based on the copula distribution of a $d$-dimensional random vector. A $d$-dimensional copula is a probability distribution $C$ on the $d$-dimensional unit hypercube $[0,1]^d$ such that all of its one-dimensional marginals are uniform on $[0,1]$. If $F$ is the distribution function of  a $d$-dimensional random vector $\Xvec=(X_1,\ldots,X_d)$ with continuous one-dimensional marginals $F_1,\cdots, F_d$, then the copula transformation of $F$ or the copula distribution of $\Xvec$ is given by
\vspace{-.2ex}
\begin{align}
\cx(\uvec)=F(F_1^{-1}(u_1),F_2^{-1}(u_2),\cdots,F_d^{-1}(u_d)),
\label{eq:0}
\end{align}
where $\uvec=(u_1,u_2,\ldots,u_d) \in [0,1]^{d}$ and $F_i^{-1}(u_i)\!=\!\inf\{x\!:\!F_i(x)\!>\!u_i\}$ for all $i=1,2,\ldots,d$ \citep[see, e.g.,][for further discussion on copula]{nelsen2006}. If $\cx$ is the cumulative distribution function of a uniform distribution on $[0,1]^{d}$, i.e., $X_1,\ldots,X_d$ are independent, the associated copula is called the uniform copula and denoted by $\Pi$. On the other hand, if $X_1,\ldots,X_d$ are comonotonic, i.e. there exist strictly increasing functions $g_i$'s and a random variable $V$ such that $\Xvec\stackrel{d}{=}(g_1(V),g_2(V),\ldots,g_d(V))$, it is called the maximum copula and denoted by ${\mathsf M}$. So, for every $\uvec \in [0,1]^d$, we have $\Pi(\uvec) = \prod_{i=1}^{d}u_i$ and ${\mathsf M}(\uvec)=\min\{u_1,\ldots,u_d\}$.

Naturally, larger difference between $\cx$ and $\Pi$ indicates higher degree of dependence among  $X_1,\ldots,X_d$. To measure the difference between two probability distributions $\mathbb{P}$ and $\mathbb{Q}$ on $\mathbb{R}^d$, we use 
\vspace{-.2ex}
\begin{align}
\gamma_k(\mathbb{P},\mathbb{Q})=\Bigl[\E\left\{k(\Yvec,\Yvec^{'})+k(\Zvec,\Zvec^{'})-2k(\Yvec,\Zvec)\right\}\Bigr]^{1/2};
\label{eq: gam exp}
\end{align}
where $\Yvec,\Yvec^{'}\stackrel{i.i.d.}{\sim} {\mathbb P}$,  $\Zvec,\Zvec^{'}\stackrel{i.i.d.}{\sim} {\mathbb Q}$ are  independent, and $k:\Rd\times\Rd$ $\rightarrow\mathbb{R}$ is a symmetric, bounded, positive definite kernel. It is known that $\gamma_k$ is a pseudo-metric on the space of all continuous probability distributions on $\Rd$, and it is a metric when $k$ is a characteristic kernel \citep[see][]{fukumizu2008kernel}). Gaussian kernel $\ks(x,y)=\exp\left(-\frac{\|x - y\|^2}{2\sigma^2}\right)$ with a bandwidth parameter $\sigma>0$ is a popular choice as a characteristic kernel, and we shall use it throughout this article.

From the above discussion, it is clear that for any characteristic kernel $k$ on $\Rd\times\Rd$, one can use $\gamma_k({\mathsf C}_{\Xvec},\Pi)$ or  $\gamma^2_k({\mathsf C}_{\Xvec},\Pi)$ as a measure of dependency among the coordinate variables $X_1,X_2,\ldots,X_d$. In this article, we use a scaled version of this measure given by
\vspace{-.2ex}
\begin{align}
\cb={{\gamma_{\ks}(\cx,\Pi)}}/{{\gamma_{\ks}(\mathsf{M},\Pi)}}.
\label{eq:1}
\end{align}
Note that  the denominator $\gamma_{\ks}(\mathsf{M},\Pi)$ is strictly positive. So, $\cb$ is well defined. The use of the Gaussian kernel makes the measure $\cb$ invariant under permutations and strictly monotone transformations of the coordinate variables. The result is stated below.

\begin{lemma}
	\label{lem: per mon inv}
	$\cb$ is invariant under permutations and strictly monotone transformations of $X_1,X_2,\ldots,X_d$.
\end{lemma}	

From the definition of $\cb$, it is clear that it takes the value $0$ if and only if the coordinates of $\Xvec$ are independent, and its value is supposed to increase as the dependency among $X_1,X_2,\ldots,X_d$ increases. The following lemma shows that in case of extreme dependency (i.e., when for each pair of variables, one is a strictly monotone function of the other), it turns out to be $1$.

\begin{lemma}
	\label{lem: high val}
	For all $i=1,\ldots,d-1$, if $X_i$ is a strictly monotone function of $X_d$ with probability one, then  $\cb$  takes the value $1$.
\end{lemma}

This desirable property of $\cb$ helps us to properly assess the degree of dependency among $X_1,\ldots,X_d$.
Note that many well-known dependency measures like the copula based multivariate extensions of Spearman's $\rho$, Kendall's $\tau$, Blomqvist's $\beta$ and Hoeffding's $\phi$ statistic \citep[see, e.g.,][]{flores05,nelsen1996nonparametric,gaisser2010multivariate} do not have this property. 
But, since $\cb$ always takes value in $[0,1]$, in the case of $d=2$, unlike Spearman's $\rho$, Kendall's $\tau$ statistic or  Blomqvist's $\beta$ statistic, it does not give us any idea about the direction of dependence between two variables. 
We know that the distance correlation measure proposed by \cite{szekely2007measuring} can be expressed as a weighted squared distance between the characteristic functions of two distributions. The following theorem shows that 
$\cb$ also has a similar property. 
\begin{theorem}
	\label{thm: rep}
	Let $\varphi_{\cx}$ and $\varphi_{\Pi}$ be the characteristic functions of $\cx$ and $\Pi$, respectively.  Define $\csd=\kappa\left(\frac{\sigma}{\sqrt{d}}\right)+\kappa^d(\sigma)-2\int\displaylimits_0^1 \lambda^d(u,\sigma)\,du,$ where\\ $\kappa(\sigma)=\sqrt{2\pi}\sigma\left[2\Phi\left( \frac{1}{\sigma}\right)-1\right]-2\sigma^2\left[1-\exp\left(-\frac{1}{2\sigma^2}\right) \right]$,\\ $\lambda(x,\sigma)=\sqrt{2\pi}\sigma\left[  \Phi\left(\frac{x}{\sigma}\right)+\Phi\left(\frac{1-x}{\sigma}\right)-1\right]$ and $\Phi(\cdot)$ is the cumulative distribution function of the N(0,1)  distribution.
	Then $\scb$ can be expressed as
	$$\scb=\csd^{-1}\left(\frac{\sigma}{\sqrt{2\pi}}\right)^d\int_{\mathbb{R}^d}|\varphi_{\cx}(\wvec)-\varphi_{\Pi}(\wvec)|^2\exp\left(-\frac{\sigma^2}{2}\wvec^{\intercal}\wvec\right)\,d\wvec.$$
\end{theorem}

Another interesting property of $\cb$ is its continuity. Note that if $\{\Xvec_n; n\ge 1\}$ is a sequence of random vectors converging in distribution to $\Xvec$, then $\mathsf{C_{\Xvec_n}}$ converges to $\cx$ weakly. So, using the dominated convergence theorem, from Theorem 1 it follows that $I_\sigma(\Xvec_n)$ converges to $\cb$ as $n$ increases. This result is stated below.

\begin{lemma}
	\label{lem: cont}
	Let $\{\Xvec_n: n\ge 1\}$ be a sequence of $d$-dimensional random vectors with continuous one-dimensional marginals. If $\Xvec_n$ converges to $\Xvec$ weakly, we have $\lim_{n\to\infty}I_\sigma(\Xvec_n)=\cb.$
\end{lemma}

In the case of $d=2$, $\cb$ enjoys some additional properties. For instance, $\scb$ can be viewed as a product moment correlation coefficient between two random quantities. 
 If $\Xvec$ follows a bivariate normal distribution with correlation coefficient $r$, $\cb$ turns out to be a strictly increasing function of $|r|$. These results are shown by the following theorem.

\begin{theorem}
	\label{thm: d2 rep}
	Let $\Xvec=(X_1,X_2)$ be a bivariate random vector with continuous one-dimensional marginals.
	\begin{enumerate}
		
		\item[$(a)$] Suppose that $\Tvec=(T_1,T_2)$ and $\Tvec^{'}=(T_1^{'},T_2^{'})$ are  independent, and they follow the distribution $\cx$, the copula distribution of $\Xvec$. Define $V_i=\ks(T_i,T_i^{'})-\E\left[\ks(T_i,T_i^{'})\Big|T_i\right]-\E\left[\ks(T_i,T_i^{'})\Big|T_i^{'}\right]+\E\left[ \ks(T_i,T_i^{'})\right]$ for $i=1,2$. Then we have $\scb=\cor(V_1,V_2)$, which takes the value $1$ if only if $X_1$ is a strictly monotone function of $X_2$ with probability one.
		\item[$(b)$] If	$\Xvec$ follows a bivariate normal distribution with $\cor(X_1,X_2)=r$, then $\cb$ is a strictly increasing function of $|r|$ with $\cb\leq|r|.$
	\end{enumerate}
\end{theorem}

Another interesting property of $\cb$ is its irreducibility. Following \cite{Schmid2010copula}, we call a dependency measure $I$ irreducible if, for any $d>2$, $I(X_1,\ldots , X_d)$ is not a function of the quantities $\{I(X_{i_1},\ldots, X_{i_k}): \{ i_1,\cdots , i_k\}\subsetneqq\{1,\ldots ,d\}\}$. Naturally, any reasonable multivariate measure of dependency is expected to be irreducible. Note that if $I(X_1,X_2,X_3)$ gets completely determined by $I(X_1,X_2)$, $I(X_2,X_3)$ and $I(X_3,X_1)$, instead of mutual dependence among $X_1$, $X_2$ and $X_3$, it can only detect pairwise dependence. 
The following theorem shows that any copula based multivariate dependency measure, which takes the value zero only for the uniform copula, is irreducible.

\begin{theorem}
	\label{thm: irr}
	Let $\cx$ be the copula distribution of $\Xvec$ and $I (\Xvec)=\mathcal{M}(C_{\Xvec})$ be a copula based multivariate dependency measure. If $\mathcal{M}(\cx)=0$ implies $\cx=\Pi$,  then $I$ is irreducible.
\end{theorem}

For any fixed bandwidth parameter $\sigma$, the irreducibility of our proposed measure $I_{\sigma}(\Xvec)$ follows from Theorem \ref{thm: irr} as a corollary. However, this property vanishes when $\sigma$ diverges to infinity. In such a situation, the limiting value of $I_{\sigma}(\Xvec)$ turns out to be the average of squared Spearman's rank correlations between $\binom{d}{2}$ pairs of random variables as stated in the following theorem.

\begin{lemma}
	\label{lem: inf band}
	As $\sigma$ diverges to infinity, $\scb$ converges to $\frac{1}{{d \choose 2}}\sum_{1\leq i< j\leq d}\cor^2(S_i,S_j)$, where $\Svec=(S_1,S_2,\cdots,S_d)\sim \cx$.
\end{lemma}

\section{Estimation of the proposed measure}

Let $\Xvec^{(1)}, \cdots, \Xvec^{(n)}$ be $n$ independent copies of the random vector $\Xvec$ taking values in ${\mathbb R}^d$. For any fixed $j=1,\ldots,d$, define $R_j^{(i)}$ as the rank of $X_j^{(i)}$ ($i=1,\ldots,n$) in the set $\{X_j^{(1)},\ldots,X_j^{(n)}\}$ to get $\mathbf{R}^{(i)}=(R_1^{(i)},\ldots,R_d^{(i)})$, the coordinate-wise rank of $\Xvec^{(i)}$. We use the normalized rank vectors $\Yvec^{(i)}=\Rvec^{(i)}/n~$  ($i=1,\ldots,n$) to define the empirical version of the copula distribution $\cx$, which is given by
\begin{align}
&\cxn(u_1,u_2,\ldots,u_d)=\frac{1}{n}\sum_{i=1}^{n}{\mathbb I}[Y_1^{(i)}\leq u_1,Y_2^{(i)}\leq u_2,\ldots,Y_d^{(i)}\leq u_d],
\label{eq: 3}
\end{align}
where ${\mathbb I}$ is the indicator function. 	Clearly, $\cxn$ is the empirical distribution function based on $\Yvec^{(1)},\cdots,\Yvec^{(n)}$. Indeed, it is the copula transform of the empirical distribution based on $\Xvec^{(1)},\ldots,\Xvec^{(n)}$. Similarly, we can define empirical versions of the maximum copula and the uniform copula as
\begin{align}
&\mathsf{M}_n(u_1,u_2,\ldots,u_d)=\frac{1}{n}\sum_{i=1}^{n}{\mathbb I}[u_1 \geq i/n,u_2 \geq i/n,\ldots,u_d \geq i/n] ~\mbox{and}\\ \nonumber
&\Pi_n(u_1,u_2,\cdots,u_d)=\frac{1}{n^d}\sum_{1\leq i_1,i_2,\cdots,i_d\leq n}{\mathbb I}[u_1 \geq i_1/n,u_2 \geq i_2/n,\cdots,u_d \geq i_d/n],
\label{eq: 4}
\end{align}
respectively. While $\mathsf{M}_n$ puts mass $1/n$ on each of the $n$ points $\{({i}/{n},{i}/{n},\cdots,{i}/{n}):$ $~1\leq i\leq n\}$, $\Pi_n$ assigns equal mass to $n^d$ points of the form $({i_1}/{n},{i_2}/{n},\cdots,{i_d}/{n})$, for $i_1,i_2,\cdots,i_d\in\{1,2,\cdots,n\}$. We estimate   $\cb$ by its empirical analog

\begin{equation}
\label{eq: 5}
\ecb={\gamma_{\ks}(\cxn,\Pi_n)}/{\gamma_{\ks}(\mathsf{M}_n,\Pi_n)}.
\end{equation}

Note that $\ecb$ is well-defined since $\mathsf{M}_n\not=\Pi_n$ for every $n>1$. One can also check that $\ecb$ can be expressed as (replace expectation by sample mean in equation (2)) 
\begin{equation}
\label{eq: ecb form}
\ecb=\sqrt{\dfrac{s_1-2s_2+v_3}{v_1-2v_2+v_3}},
\end{equation}
where $s_1=\frac{1}{n^2}\sum\limits_{1\leq i<j\leq n}^{} k_{\sigma}(\Yvec^{(i)},\Yvec^{(j)})+\frac{1}{n},\; s_2=\frac{1}{n^{d+1}}\sum\limits_{i=1}^{n}\prod\limits_{j=1}^{d}\sum\limits_{l=1}^{n}\ks(Y_j^{(i)},\frac{1}{n}),$ $v_1=\frac{2}{n^2}\sum\limits_{i=1}^{n-1}(n-i)e^{-\frac{d}{2}\left({i}/{n\sigma}\right)^2}+\frac{1}{n},\; v_2=\frac{1}{n^{d+1}}\sum\limits_{i=1}^{n}\left[\sum\limits_{j=1}^{n}e^{-\frac{1}{2}\left({(i-j)}/{n\sigma}\right)^2}\right]^d$ and  $v_3=\left[\frac{2}{n^2}\sum\limits_{i=1}^{n-1}(n-i)e^{-\frac{1}{2}\left({i}/{n\sigma}\right)^2}+\frac{1}{n}\right]^d$.\\

The above formula shows that the computing cost of $\ecb$ is $O(dn^2)$. This estimate enjoys some nice theoretical properties similar to those of $I_{\sigma}(\Xvec)$. These properties are mentioned below.

\begin{lemma}
	\label{lem: emp prop}
	Suppose that $\Xvec^{(1)},\Xvec^{(2)},\ldots,\Xvec^{(n)}$ are $n$ independent observations from the distribution of a $d$-dimensional random vector $\Xvec=(X_{1},X_{2},\ldots,X_{d})$ having continuous one-dimensional marginals. Then, we have the following results.
	\begin{enumerate}
		\item[(a)] $\ecb$ is invariant under permutation and strictly monotone transformations of the coordinate variables  $X_{1},X_{2},\ldots,X_{d}$.
		\item[(b)]  For all $i=1,\ldots,d-1$, if $X_i$ is a strictly monotone function of $X_d$ with probability one, then $\ecb$ takes the value $1$.
	\end{enumerate}
\end{lemma}

Note that other existing copula based dependency measures
do not have the property mentioned in part ($b$) of Lemma 5.
For instance, multivariate extensions of Spearman's $\rho$,  Kendall's $\tau$, Blomqvist's $\beta$ and Hoeffding's $\phi$ statistics
\citep[see, e.g.,][]{nelsen1996nonparametric,nelsen2002concordance,gaisser2010multivariate} may not take the value $1$ even when the measurement variables have monotone relationships among them. To demonstrate this, we considered a simple example. We generated 10000 observations on $\Xvec=(\Xd)$, where $X_i=V$ or $X_i=-V$ for $V\sim U(0,1)$ and $i=1,2,\ldots,d$. Hence each pair of variables were monotonically related. We considered three choices of $d$ ($d=3,4,5$), and for each value of $d$, results are reported in Table \ref{tab: mono} for different types of relationships shown in the orientation column. For example, the $(\uparrow,\uparrow,\downarrow)$ sign in the orientation column indicates that $(X_1,X_2,X_3)=(V,V,-V)$. Table \ref{tab: mono} clearly shows that all dependency measures considered here take the value $1$ when the relationships among the variables are strictly increasing. But, except for $\ecb$, all other measures fail to have this property for other monotone relationships among the variables.

\begin{table}[ht]
	\begin{center}
		\begin{tabular}{|c c |c c c c c|}
			\hline
			Dimension & Orientation & $\ecb$ & $\widehat{\rho}$ & $\widehat{\tau}$ & $\widehat{\beta}$ & $\widehat{\varphi}$ \\
			\hline
			3 & $(\uparrow,\uparrow,\uparrow)$ & 1.000 & 1.000 & 1.000 & 1.000 & 1.000\\
			3 & $(\uparrow,\uparrow,\downarrow)$ & 1.000 & -0.333 & -0.333 & -0.333 & 0.517\\
			4 & $(\uparrow,\uparrow,\uparrow,\uparrow)$ & 1.000 & 1.000 & 1.000 & 1.000 & 1.000\\
			4 & $(\uparrow,\uparrow,\uparrow,\downarrow)$ & 1.000 & -0.091 & -0.143 & -0.143 & 0.382\\
			4 & $(\uparrow,\uparrow,\downarrow,\downarrow)$ & 1.000 & -0.212 & -0.143 & -0.143 & 0.327\\
			5 & $(\uparrow,\uparrow,\uparrow,\uparrow,\uparrow)$ & 1.000 & 1.000 & 1.000 & 1.000 & 1.000\\
			5 & $(\uparrow,\uparrow,\uparrow,\uparrow,\downarrow)$ & 1.000 & 0.016 & -0.067 & -0.067 & 0.347\\
			5 & $(\uparrow,\uparrow,\uparrow,\downarrow,\downarrow)$ & 1.000 & -0.108 & -0.067 & -0.067 & 0.236\\
			\hline
		\end{tabular}
		\caption{Different measures of dependency when the variables are monotonically related.}
		\label{tab: mono}
	\end{center}
	\vspace{-0.25in}
\end{table}

Since $\ecb$ is based on coordinate-wise ranks of the observations, it is robust against contaminations and outliers generated from
heavy-tailed distributions. Following the results in \cite{poczos2012}, one can show that addition of a new observation can change its value
by at most ${\bf O}(n^{-1})$. Just like $\cb$, its empirical analog $\ecb$ also enjoys some additional properties for $d=2$. Theorem \ref{thm: est upbnd} below shows that result analogous to Theorem \ref{thm: d2 rep}(a) holds for $\ecb$ as well. 

\begin{theorem}
	\label{thm: est upbnd}
	Suppose that $\Xvec^{(1)},\cdots,\Xvec^{(n)}$ are independent observations from a bivariate distribution with continuous one-dimensional marginals and $\Yvec^{(1)},\cdots,\Yvec^{(n)}$ are their normalized coordinate-wise ranks.
	For $i,j=1,2,\ldots,n$, define
	\[
	\begin{array}{ll}
	V_1(i,j)=\ks(Y_1^{(i)},Y_1^{(j)})-\frac{1}{n}\sum_{i=1}^{n}\ks(Y_1^{(i)},Y_1^{(j)})-\frac{1}{n}
	\sum_{j=1}^{n}\ks(Y_1^{(i)},Y_1^{(j)})\\
	\qquad\qquad+\frac{1}{n^2}\sum_{i,j=1}^{n}\ks(Y_1^{(i)},Y_1^{(j)}),
	~\mbox{and}\\
	V_2(i,j)=\ks(Y_2^{(i)},Y_2^{(j)})-\frac{1}{n}\sum_{i=1}^{n}\ks(Y_2^{(i)},Y_2^{(j)})-\frac{1}{n}
	\sum_{j=1}^{n}\ks(Y_2^{(i)},Y_2^{(j)})\\
	\qquad\qquad+\frac{1}{n^2}\sum_{i,j=1}^{n}\ks(Y_2^{(i)},Y_2^{(j)}).
	\end{array}
	\]
	Then $\ecb$ can be expressed as\\ $\ecb = \Big[\sum_{i,j} V_1(i,j)V_2(i,i)\Big]/ \Big[\sum_{i,j} V_1^2(i,j) \sum_{i,j} V_2^2(i,j)\Big]^{1/2}$.
	As a consequence, we have  $0 \leq \ecb \leq 1$, where the equality holds if and only if one coordinate variable is a  strictly monotone function of the other.
\end{theorem}

\section{Test of independence based on $\ecb$}
\label{sec: mul toi}

We have seen that $\cb$ serves as a measures of dependence among the coordinates of $\Xvec$. It is non-negative, and takes the value $0$ if and only if $X_1,X_2,\ldots,X_d$ are independent. So, we can use $\ecb$ as the test statistic and reject $H_0$, the null hypothesis of independence, for large values of $\ecb$. The large sample distribution of our test statistic is given by the following theorem.

\begin{theorem}
	\label{thm: asymp conv}
	Suppose that $\Xvec$ follows a multivariate distribution with continuous one-dimensional marginals. {{Also assume that the associated copula distribution $\cx$ has continuous partial derivatives.}}
	\begin{itemize}
		\item[(a)] $\mbox{If}~ \cx=\Pi, \mbox{ then}~~ n\secb
		\overset{\mathcal{L}}{\longrightarrow}\;\sum_{i=1}^{\infty}\lambda_iZ_i^2$, where the
		$Z_i$'s are {i.i.d.} $\mathcal{N}(0,1)$ and the $\lambda_i$'s are some positive constants (see the proof of the theorem in the Appendix for detailed description of the $\lambda_i$'s).
		
		\item[(b)] If $\cx\neq\Pi$, then
		$\sqrt{n}(\ecb-\cb)\overset{\mathcal{L}}{\longrightarrow}\mathcal{N}(0,\delta^2)$; where\\ $\delta^2\!=\!\csd^{-2}I^{-2}_\sigma(\Xvec)\!\int\limits_{[0,1]^d}\!\int\limits_{[0,1]^d}\!g(u)g(v)\,\E[\,d\mathbb{G}_{\cx}(u)\,d\mathbb{G}_{\cx}(v)]$,\\ $g(u)\!=\!\!\int\limits_{[0,1]^d}\!\!\ks(u,v)d(\cx\!-\Pi)(v)$ and $\mathbb{G}_{\cx}$ is a 0 mean Gaussian process (as defined in Theorem T1 in the Appendix).
	\end{itemize}
\end{theorem}

\begin{figure}[h]
	\vspace{-0.2in}
	\centering
	\includegraphics[width=0.9\textwidth]{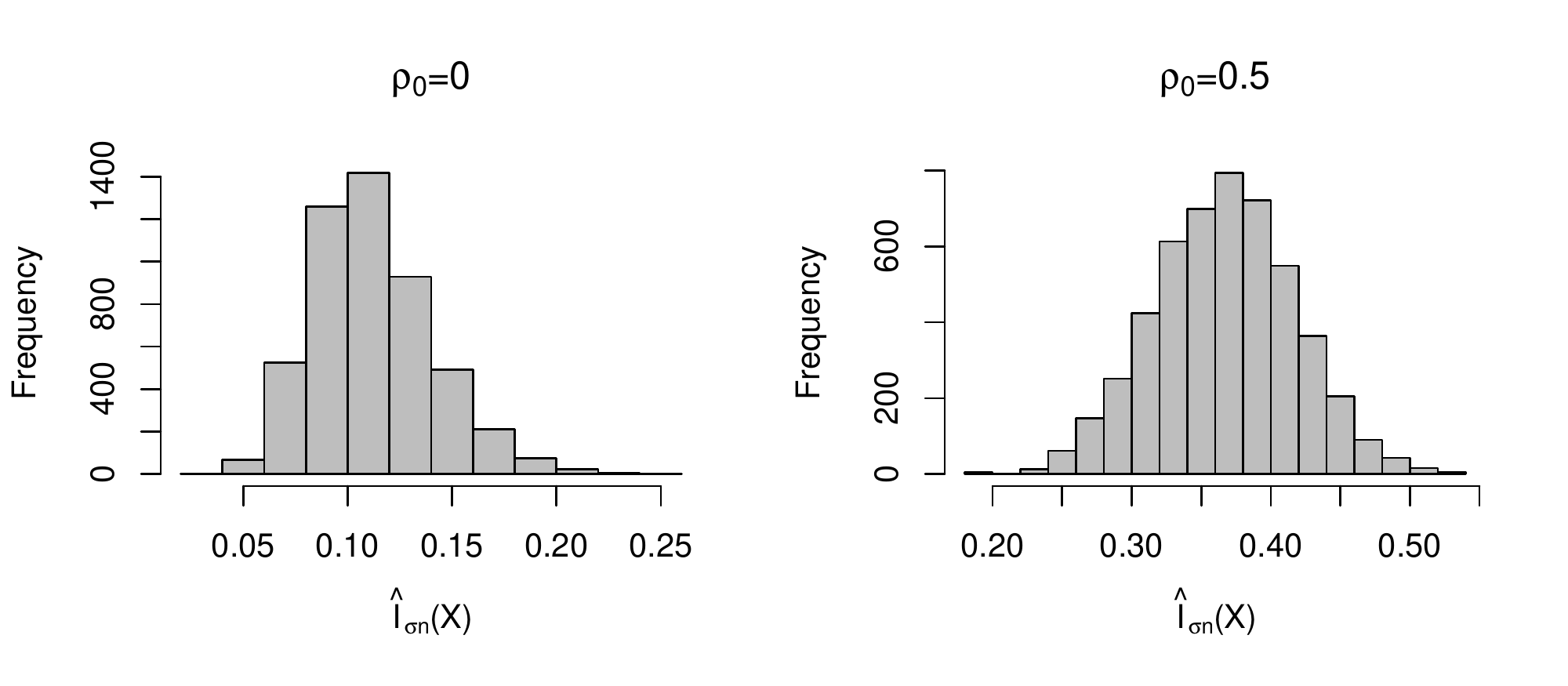}
	\vspace{-0.2in}
	\caption{Empirical distribution of $\ecb$ for $\sigma=0.2$}
	\label{fig: empdist}
	\vspace{-0.15in}
\end{figure}

The histograms in Figures \ref{fig: empdist}(a) and \ref{fig: empdist}(b) show the empirical distributions of $\ecb$ computed based on 5000 independent samples, each of size 200, generated from bivariate normal distributions with correlation coefficient $\rho_0=0$ and $0.5$, respectively. For $\rho_0=0.5$ (i.e., $\cx \neq \Pi$) while the empirical distribution looks like a normal distribution,  for $\rho_0=0$ (i.e., $\cx= \Pi$), it turns out to be positively skewed. This is consistent with the result stated in Theorem \ref{thm: asymp conv}.

The probability convergence of $\ecb$ follows from Theorem \ref{thm: asymp conv}. But, we also have a stronger result in this context. The following theorem shows that $\ecb$ converges to 
$\cb$ almost surely.

\begin{theorem}
	\label{thm: as conv}
	$\ecb$ converges to $\cb$ almost surely as the sample size $n$ tends to infinity.
\end{theorem}

From Theorem \ref{thm: as conv}, it is clear that under the null hypothesis of independence, $\ecb$ converges to $0$ almost surely, while under the alternative, it converges to a positive constant. For any fixed choice of $\sigma$, the large sample consistency of the test follows from it. 
However, for practical implementation of the test, one needs to determine the cut-off. It is difficult to find this cut-off based on the asymptotic null distribution of the test statistic mentioned in Theorem \ref{thm: asymp conv} since the coefficients $\lambda_i$'s associated with the chi-square distributions are all unknown. 
Here we use the distribution-free property of $\ecb$ to determine the cut-off. Note that under $H_0$, for each $j=1,2,\ldots,d$, we have $P(R_j^{(1)}=i_1,R_j^{(2)}=i_2,\ldots,R_j^{(n)}=i_n)=1/n!$ for any permutation $\{i_1,i_2,\ldots,i_n\}$ of $\{1,2,\ldots,n\}$, and for different values of $j$, they are independent. So, we can easily generate normalized coordinate-wise ranks to compute the test statistic. We repeat this procedure 10,000 times to approximate the $(1-\alpha)$-th quantile of the null distribution of $\ecb$, which is then used as the cut-off. This whole calculation can be done off-line, and we can prepare a table of critical values for different choices $n$ and $\sigma$ before handling actual observations.

Though any fixed choice of the bandwidth $\sigma$ leads to a consistent test (follows from Theorem \ref{thm: as conv}), its power may depend on this choice. The method commonly used for choosing the bandwidth is based on ``median heuristic" \citep[see, e.g.,] [Sec~5]{fukumizu2009}, where one computes all pairwise distances among the observations and then the median of those distances is taken as the bandwidth. Since we are using the kernel on the normalized rank vectors ${\bf Y}^{(1)},{\bf Y}^{(2)},\ldots,{\bf Y}^{(n)}$ having the null distribution $\Pi_n$, following this idea, we  can choose $\sigma$ to be the median of $\|\Zvec-\Zvec^{'}\|$, where $\Zvec,\Zvec^{'} \stackrel{i.i.d.}{\sim}\Pi_n$, and $\|\cdot\|$ denotes the usual Euclidean norm. Note that the bandwidth chosen in this way is non-random and it is a function of $n$. We denote it by $\sigma_n$. As $n$ increases, since $\Pi_n$ converges to $\Pi$, $\sigma_n$ converges to the median of $\|\Zvec-\Zvec^{'}\|$, where $\Zvec,\Zvec^{'} \stackrel{i.i.d}{\sim} \Pi$. Our test remains consistent for such choices of the bandwidth. This result is stated below. 

\begin{theorem}
	\label{thm: pow 1}
	Suppose that $\Xvec^{(1)},\Xvec^{(2)},\ldots,\Xvec^{(n)}$ are independent copies of $\Xvec$, which follows a multivariate distribution with continuous univariate marginals and $\cx \neq \Pi$. Also consider a sequence of bandwidths $\{\sigma_n:~n\ge 1\}$
	converging to some $\sigma_0>0$. Then, power of the proposed test based on ${\widehat I}_{\sigma_n,n}(X)$ converges to $1$ as the sample size $n$ diverges to infinity.	
\end{theorem}

In our experiments, we observed that median heuristic performs well when the relationships among the variables are nearly monotone. 
But in cases of complex non-monotone relationships, use of smaller bandwidths often yields better results. In such cases, instead of median, one can use lower quantiles of pairwise distances.

\begin{figure}[h]
	\vspace{-0.1in}
	\centering
	\includegraphics[width=0.85\textwidth]{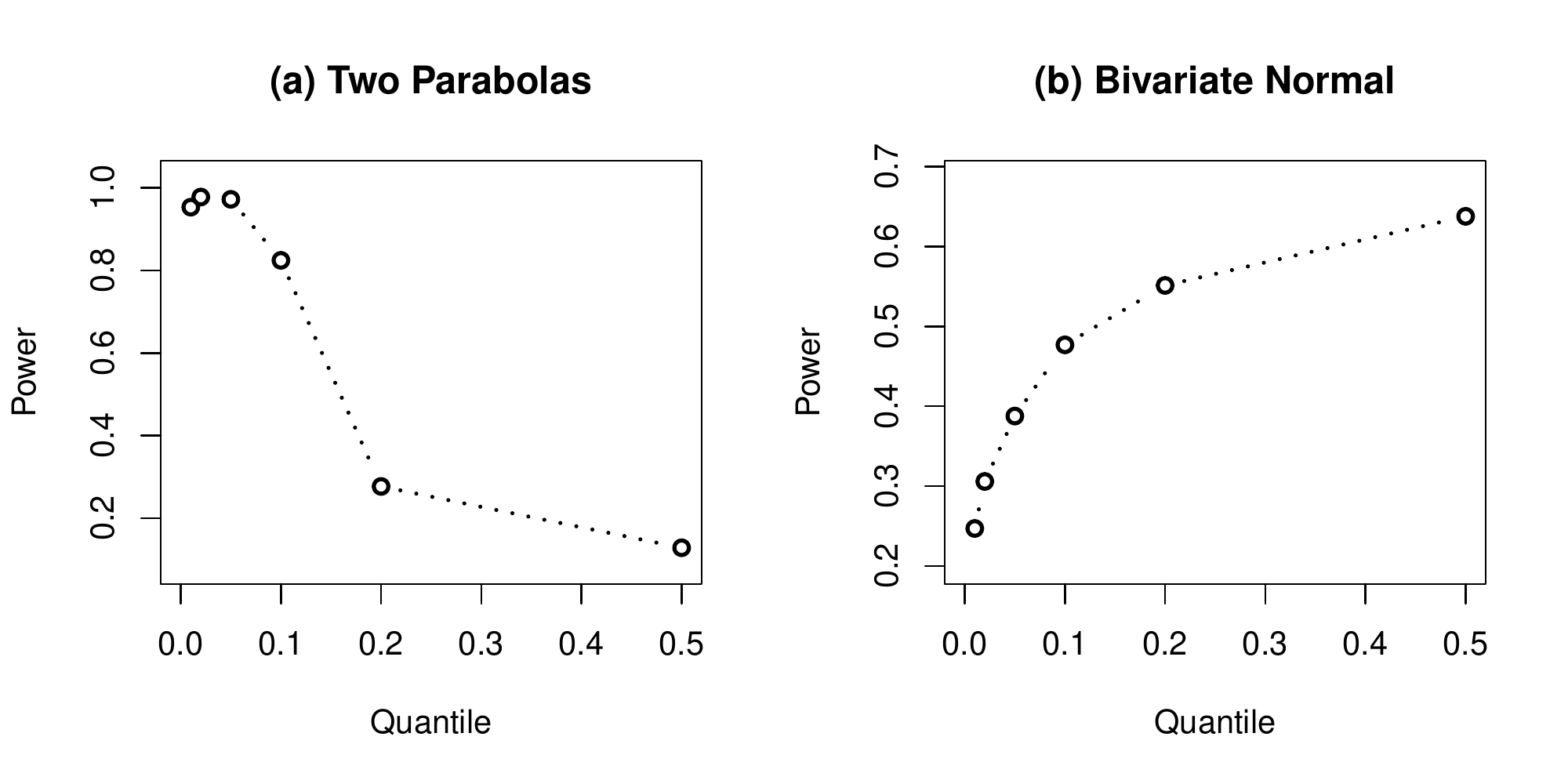}
	\vspace{-0.2in}
	\caption{Powers of the test for bandwidths based on  different quantiles of pairwise distances.}
	\label{fig: qp}
	\vspace{-0.15in}
\end{figure}

To demonstrate this, we considered two simple examples involving bivariate data sets. In one example,  observations were generated from the `Two parabolas'-type distribution mentioned in \cite{newton2009} (see Figure \ref{fig: toidd}(e)) and in the other example, they were generated from a bivariate normal distribution with correlation coefficient 0.5. In each case, we generated $25$ observations and repeated the experiment $10000$ times to estimate the powers of the tests based on $\ecb$ for different choices of $\sigma$ based on different quantiles ($0.01, 0.02, 0.05, 0.1, 0.2$ and $0.5$) of pairwise distances. Figure \ref{fig: qp} clearly shows that though median of pairwise distances worked well in the second example, smaller quantiles had better results in the first.
This figure clearly shows that depending on the underlying distribution of $\Xvec$, sometimes we need to use larger bandwidth, whereas sometimes smaller bandwidths may perform better. While larger bandwidths successfully detect global linear or monotone relationships among the variables, smaller bandwidths are useful for detecting non-monotone or local patterns. In order to capture both types of dependence, here we adopt a multi-scale approach, where we look at the results for several choices of bandwidth and then aggregate them judiciously to come up with the final decision. 

\begin{figure}[h]
	\vspace{-0.2in}	
	\centering
	\includegraphics[width=0.85\textwidth]{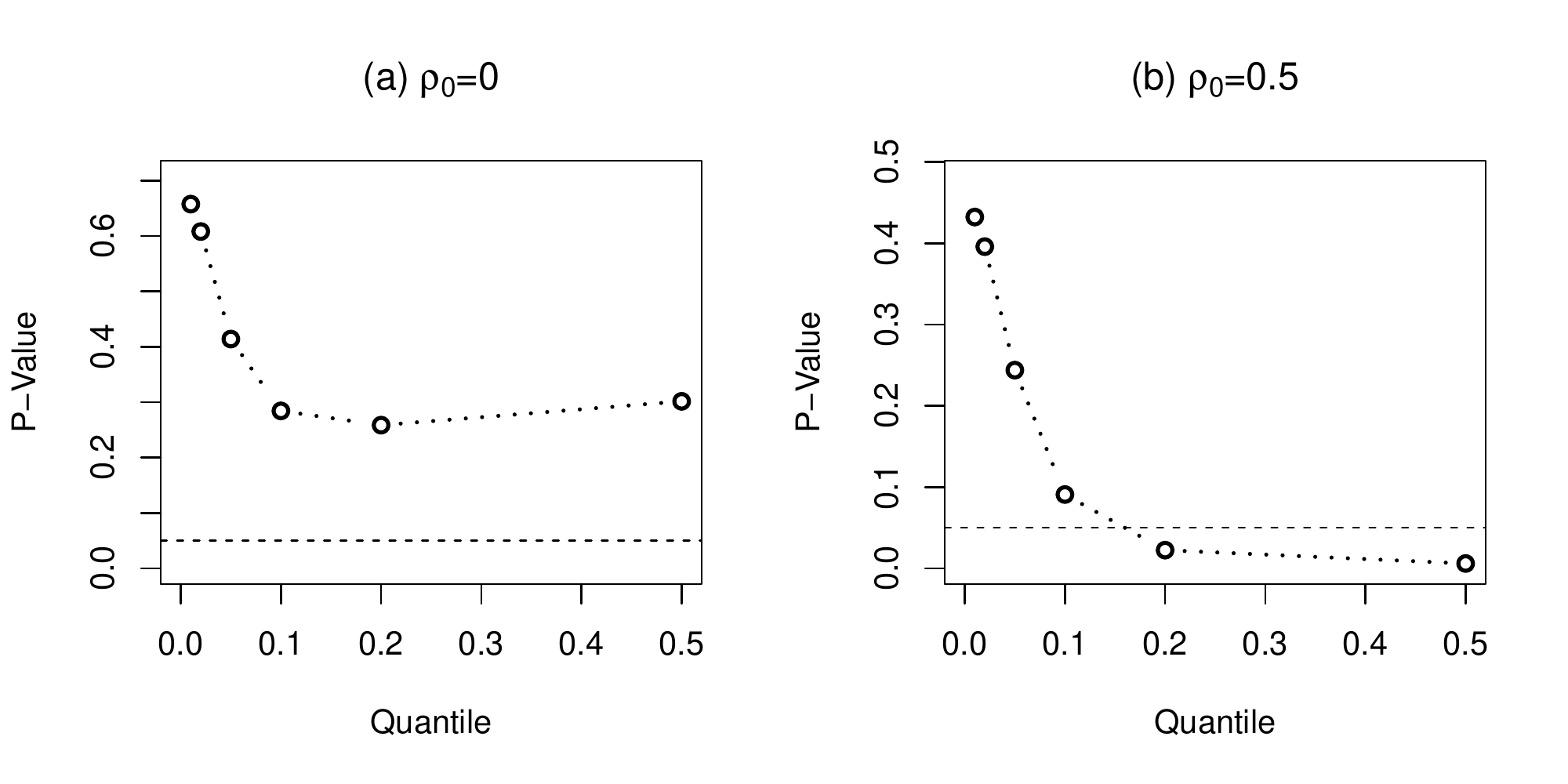}
	\vspace{-0.2in}
	\caption{p-values for bandwidths corresponding to different quantiles of pairwise distance.}
	\label{fig: pp}
	\vspace{-0.15in}	
\end{figure}

Figures \ref{fig: pp}(a) and \ref{fig: pp}(b) show the observed p-values for different choices of the bandwidth (based on quantiles of pairwise distances) when a sample of size 25 was generated from bivariate normal distributions with correlation coefficient 0 and 0.5, respectively.  Clearly, these plots of p-values carry more information than just the final result. In the first case, higher $p$-values for all choices of the bandwidth give a visual evidence in favor $H_0$, while smaller $p$-values for a long range of bandwidths in the second case indicates dependence between the two coordinate variables. Also the pattern of p-values can reveal the structure of dependence among the variables. For instance, smaller p-values for larger bandwidths indicate that the relationship between the two variables is nearly monotone, while those for smaller bandwidths indicate complex, non-monotone relations.

One way of aggregating the results corresponding to $m$ different bandwidths $\sigma^{(1)},\ldots,\sigma^{(m)}$
is to use $T_{Sum}=\sum_{i=1}^{m} {\widehat I}_{\sigma^{(i)},n}(\Xvec)$ or $T_{Max}=\max_{1\le i \le m}{\widehat I}_{\sigma^{(i)},n}(\Xvec)$ as the test statistic. Following \cite{sarkar2018some}, one can also use another method based on false discovery rate (FDR). Let $p_i$ be the $p$-value of the test based on $\sigma^{(i)}$  ($i=1,2,\ldots,m$) and $p_{(1)} \le p_{(2)} \le \ldots \le p_{(m)}$ be the corresponding order statistics. We reject $H_0$ at level $\alpha$ if and only if the set $\{i : p_{(i)} < i ~\alpha/m \}$ is non-empty. \cite{benjamini1995controlling} proposed this method for controlling FDR for a set $m$ independent tests. Later, \cite{benjamini2001control} showed that it also controls FDR when the tests statistics are positively regression dependent. Since we are testing the same hypothesis for different choices of the bandwidth, this method controls the level of the test as well \citep[see, e.g.,][]{cuesta2010simple}. It is difficult to prove positive regression dependence among the test statistics corresponding to different choices of bandwidth. However, all pairwise correlations (computed over 10000 simulations) among these test statistics were found to be positive in all of our numerical experiments. This gives an indication of positive regression dependence among the test statistics and thereby provides an empirical justification for using the above method. The following theorem shows the large sample consistency of the multi-scale versions of our tests based on  $T_{Sum}$, $T_{Max}$ and $FDR$.  

\begin{theorem}
	\label{thm: max sum fdr}
	Suppose that $\Xvec^{(1)},\Xvec^{(2)},\ldots,\Xvec^{(n)}$ are independent copies of $\Xvec$  following a multivariate distribution having continuous univariate marginals and $\cx \neq \Pi$. Then, the powers of the multi-scale versions of the proposed test based on $T_{Sum}$, $T_{Max}$ and FDR  converge to $1$ as the sample size tends to infinity.
\end{theorem}

\section{Results from the analysis of simulated and real data sets}
\label{sec: num ill}

We analyzed several simulated and real data sets to compare the performance of our proposed tests with some popular tests available in the literature. In particular, we considered the dHSIC test \citep{pfister2016kernel}, the mdCov test \citep{fan2017}  and the tests based on multivariate extensions of Hoeffding's $\phi$ \citep{gaisser2010multivariate} and  Spearman's $\rho$ \citep{nelsen1996nonparametric} statistics for comparison.  For the implementation of the dHSIC test, we used the R package ``dHSIC" \citep{dHSICR}, where we used the Gaussian kernel with the default bandwidth chosen based on median heuristic. For the mdCov test, we used the codes provided by the authors. Following their suggestion \citep[see][p. 198]{fan2017}, we standardized the data set and used unit bandwidth for all experiments. We used different options for the weight function available in the codes and reported the best result. For the tests based on Hoeffding's $\phi$ and Spearman's $\rho$ statistics (henceforth referred to as  the Hoeffding test and the Spearmen test, respectively), we used our own codes. For all these methods, conditional tests based on 10000 random permutations were used. We also considered the tests proposed in \cite{poczos2012}, where they suggested to compute the cut-offs based on probability inequalities. But this choice of cut-off makes the resulting tests very conservative. As a result, they had much lower powers compared to all other tests considered here. So, we decided not to report those results in this article. For our proposed tests, we started with the bandwidth based on median heuristic ($\sigma_{0.5}$, say) and considered other bandwidths of the form $(0.5)^i\times\sigma_{0.5}$ for $i=1,2,\ldots,m$, where $m=\lceil \log_2(\sigma_{0.5}/\sigma_{0.01})\rceil$, for $\sigma_{0.01}$ being the bandwidth based on the first quantile. Results for these bandwidths were aggregated using the three methods discussed in Section 4. However, overall performance of the tests based on $T_{Max}$ and FDR was superior than the test based on $T_{Sum}$. 
So, here we report the results for the tests based on $T_{Max}$ and FDR only.
Throughout this article, all tests are considered to have $5\%$ nominal level.

\subsection{Analysis of simulated data sets}

\begin{figure}[h]
	\vspace{-0.1in}
	\centering
	\includegraphics[width=0.9\textwidth]{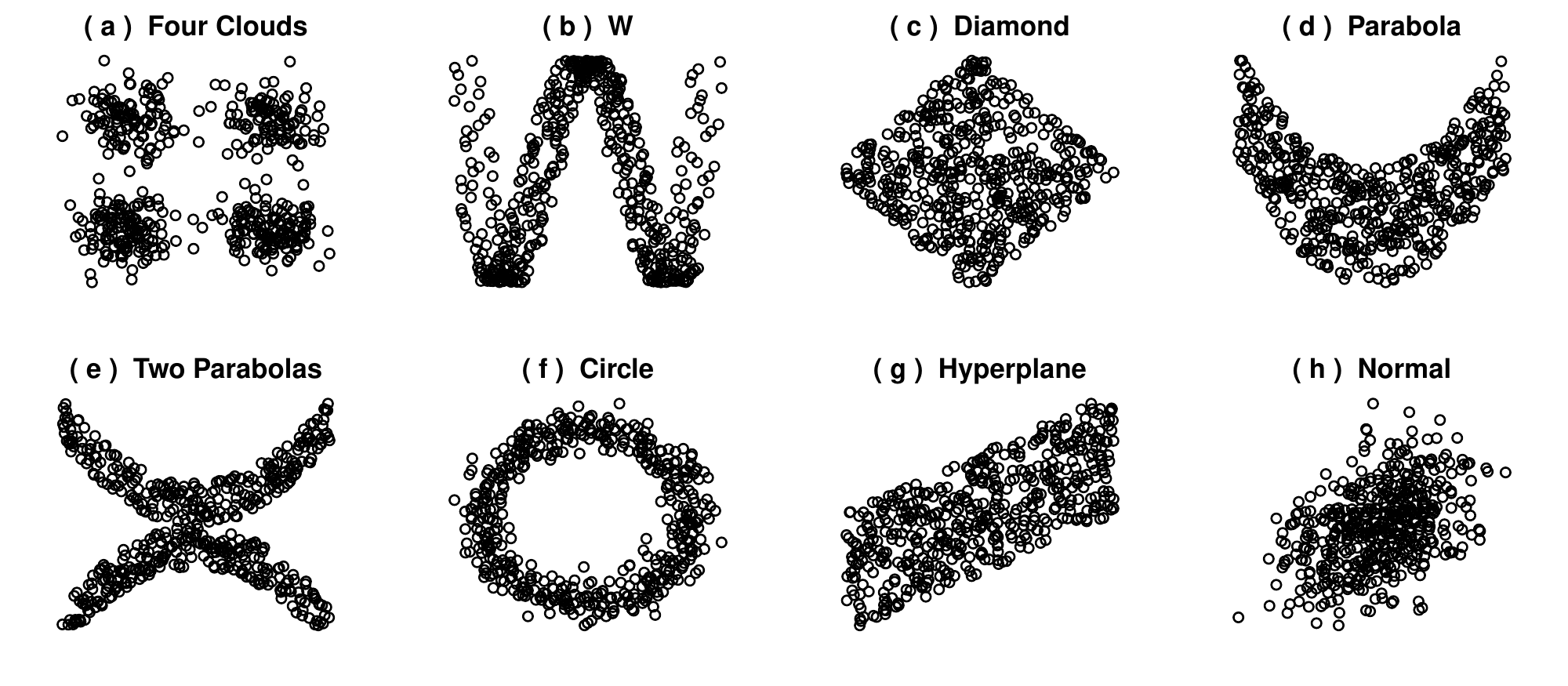}
	\vspace{-0.2in}
	\caption{Different bivariate distributions used in  simulation study.}
	\label{fig: toidd}
	\vspace{-0.15in}
\end{figure}

\begin{figure}
	\centering
	\includegraphics[width=\textwidth]{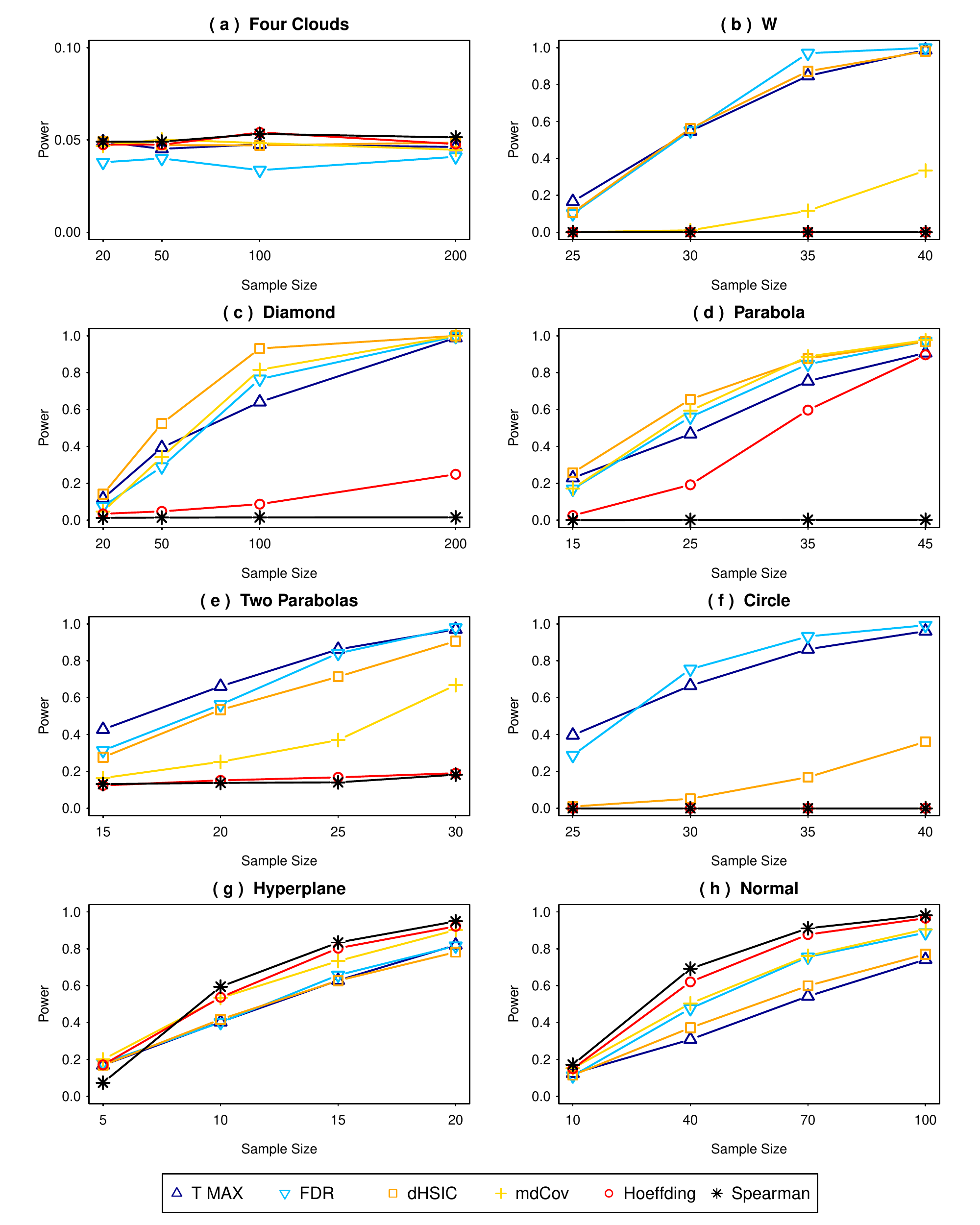}
	\vspace{-0.2in}
	\caption{Powers of different tests in simulated  bivariate data sets.}
	\label{fig: toid2}
	\vspace{-0.1in}
\end{figure}

We begin with eight simulated examples involving bivariate observations. Scatter plots of these data sets are displayed in Figure \ref{fig: toidd}. For each example, we repeated our experiment 10000 times, and the power of a test was estimated by proportion of times it rejected $H_0$. These estimated powers of different tests are reported in Figure \ref{fig: toid2}. 
The first six examples (see Figures \ref{fig: toidd}(a)-\ref{fig: toidd}(f)) are taken from \cite{newton2009}, who considered six unusual bivariate distributions. In all these examples, $X_1$ and $X_2$ are uncorrelated. In `four clouds' data, they are independent as well. In this example, almost all tests had powers close to the nomination level of 0.05 (see Figure \ref{fig: toid2}(a)). Only the test based on FDR had slightly low powers, which is quite expected in view of the conservative nature of the tests based on FDR.

In the next five examples, $X_1$ and $X_2$ are not independent. In the example with `W' type data, our proposed test based on FDR had the best overall performance followed by the test based on $T_{Max}$ and the dHSIC test (see Figure \ref{fig: toid2}(b)). Powers of all
other tests were much lower. Spearman and Hoeffding tests couldn't reject $H_0$ even on a single occasion. These two tests had zero power in `Circle'-type data as well (see Figure \ref{fig: toid2}(f)). In that example, the dCov test also had zero power, and the performance of the dHSIC test was not satisfactory as well. But our proposed test based on $T_{Max}$ and FDR performed well. These two tests outperformed their competitors in `Two parabolas'-type data as well (see Figure \ref{fig: toid2}(e)). In that example, Spearman and  Hoeffding tests again performed poorly, but the performance of mdCov and dHSIC tests was somewhat better. In the `Parabola'-type data, the test based on FDR, the dHSIC test and the mdCov test had higher powers than all other tests considered here (see Figure \ref{fig: toid2}(d)). Among the rest, the test based on $T_{Max}$ had better performance. Only in the case of `Diamond'-type data, dHSIC and mdCov tests outperformed our proposed methods (see Figure \ref{fig: toid2}(c)). 
However, even in this example, our proposed tests 
performed well. They had much higher powers compared to Spearman and Hoeffding tests.

Unlike the previous six examples, in our last two examples, $X_1$ and $X_2$ are positively correlated. In the example with `hyperplane'-type data (see Figure 3g), we have $X_1=U$ and $X_2=U+V$. where $U,V\stackrel{i.i.d.}\sim\text{U}(-1,1)$. In the  example with `normal' data (see Figure 3h), $(X_1,X_2)$ follows a bivariate normal distribution with correlation coefficient {0.4}. In these two examples, Hoeffding and Spearman tests had the best performance closely followed by the mdCov test. Our proposed tests and the dHSIC test also had competitive performance. Among these three tests, the test based on FDR had an edge.

\begin{figure}[h]
	\centering
	\includegraphics[width=0.99\textwidth]{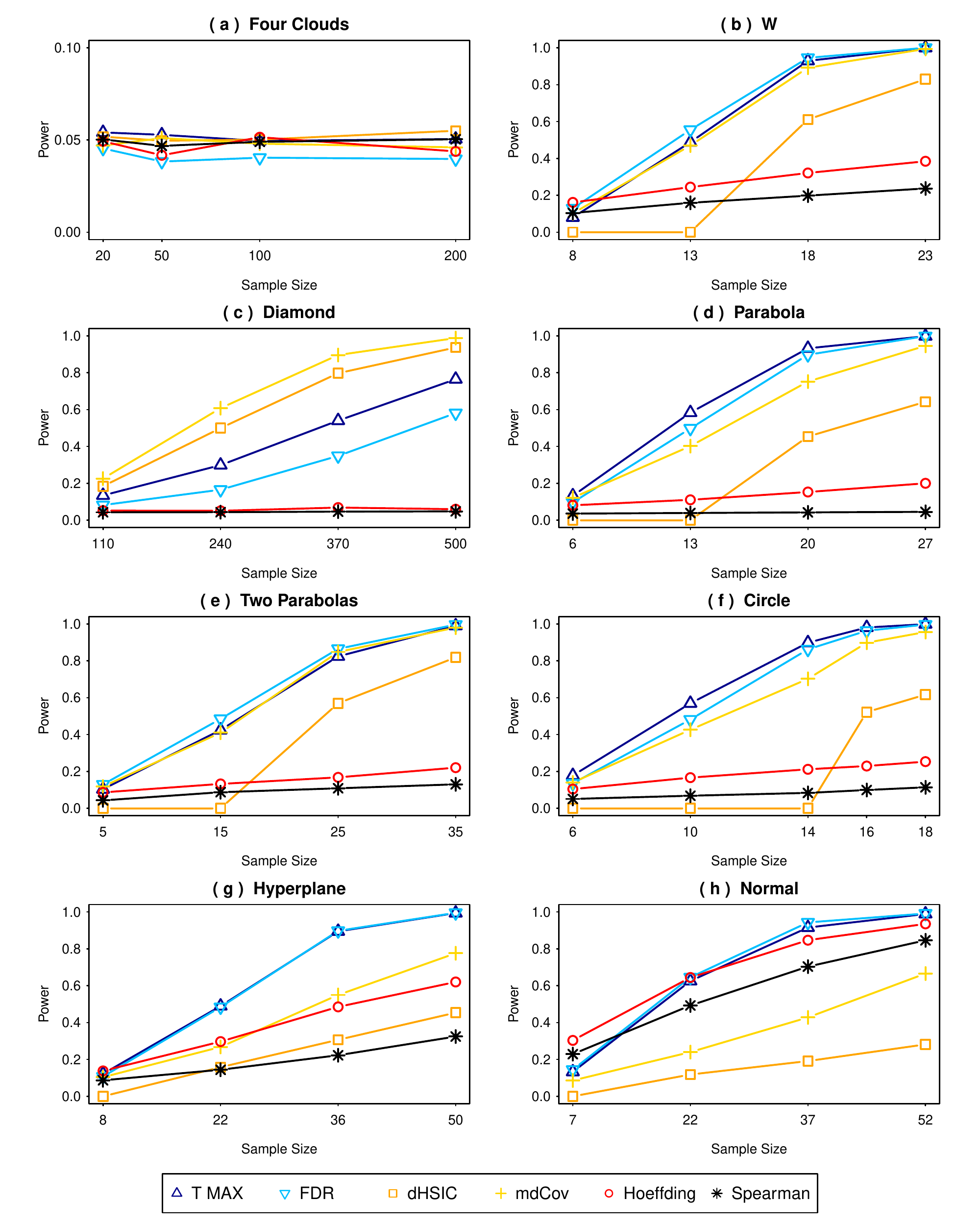}
	\vspace{-0.1in}
	\caption{Powers of different tests in eight-dimensional simulated data sets.}
	\label{fig: toid8}
	\vspace{-0.2in}
\end{figure}

Next we carried out our experiments with some eight dimensional data sets, which can be viewed  as multivariate extensions of the bivariate data sets considered above. For each of the first six examples, we generated two independent observations from the bivariate distribution, and then four independent $N(0,1)$
variables were augmented to it to get a vector of dimension eight. For the `hyperplane'-type data, we generated seven i.i.d. $N(0,1)$ variables $X_2,X_3,\cdots,X_8$, and then define $X_1=(X_2+\ldots+X_8)+\epsilon$, where $\epsilon \sim N(0,1)$. For the example with `normal' data, $\Xvec$ was generated from a $8$-dimensional normal distribution with the mean vector ${\bf 0}$ and the dispersion matrix $\mathbf{\Sigma}=((a_{i,j}))$, where $a_{i,j}=0.4^{|i-j|} ~\forall~i,j=1,2,\ldots,8$.

In the example with `four clouds' data, again the test best on FDR had powers slightly lower than 0.05, but those of all other tests were close to the nominal level (see Figure \ref{fig: toid8}(a)). 
Figure \ref{fig: toid8} clearly shows that except for `diamond'-type data, in all other cases, our tests based on $T_{Max}$ and FDR had best overall performance among the tests considered here. Note that the dHSIC test needs the sample size to be at least twice the dimension of the data (i.e., twice the number of coordinate variables) for its  implementation. So, it could not be used in some cases. In such cases, we considered its power to be zero. 

\begin{figure}
	\centering
	\includegraphics[width=0.95\textwidth]{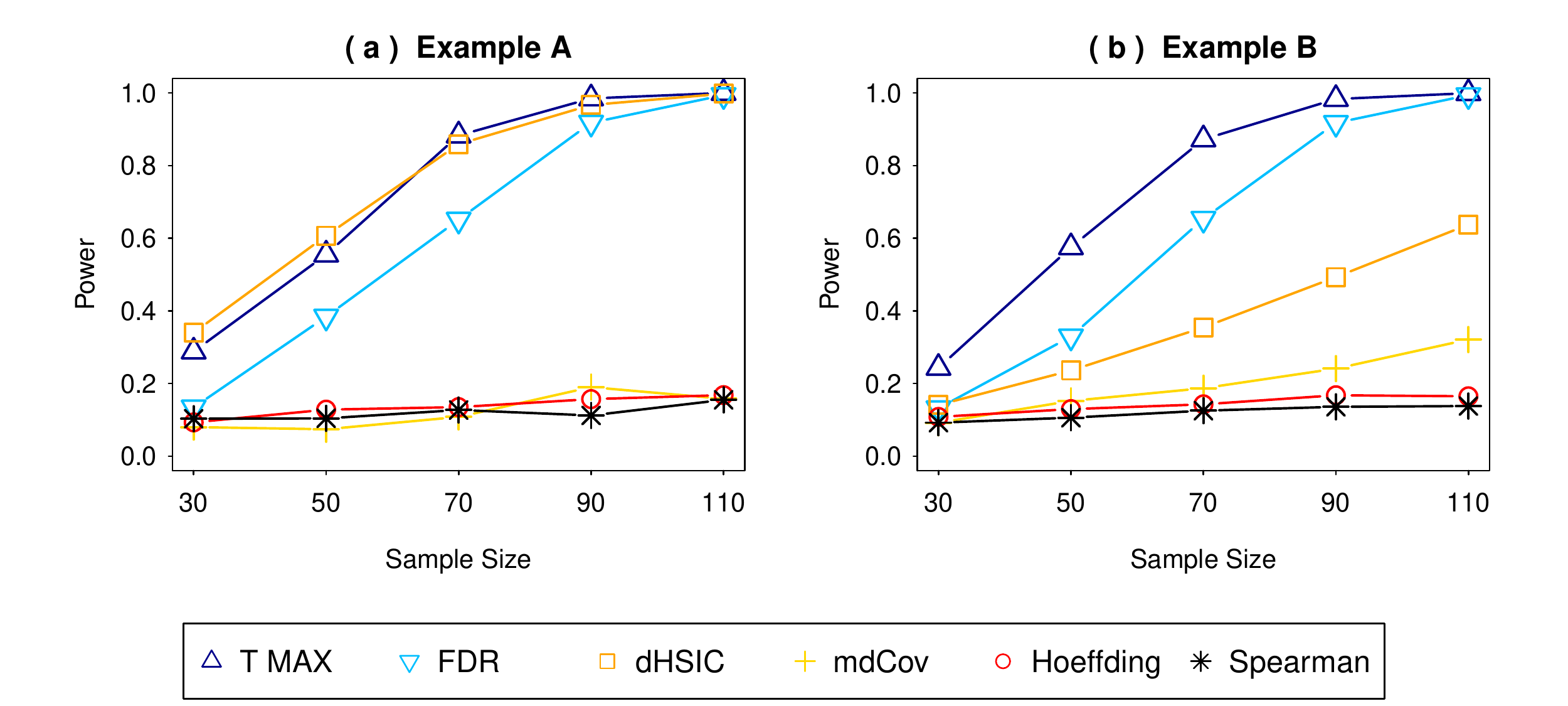}
	\caption{Powers in different tests in four-dimensional examples.}
	\label{fig: complex}
	\vspace{-0.1in}
\end{figure}

Next, we consider two interesting examples, where none of the lower dimensional marginals have dependency among the coordinate variables. In Example-A, we generate four independent $U(-1,1)$ variables $U_1,\ldots,U_4$, 
and if their product is positive, we define $X_i=U_i$ for $i=1,2,\ldots,4$. In Example-B, we generate $U_1,\ldots,U_4$ independently from $N(0,1)$ to define $X_i=U_i ~sign(U_{i+1})$ for $i=1,2,3$ and
$X_4=U_4 ~sign(U_1)$. In both of these examples, we carried out our experiments 10000 times as before to compute the powers of different tests. Note that tests based on any dependency measure, which is not irreducible, will fail to detect the dependency among the coordinate variables in these examples. Our proposed methods, particularly the test based on $T_{Max}$ had excellent performance in these two data sets (see Figure \ref{fig: complex}). In Example-A, the dHSIC test also had competitive powers, but performances of all other tests were much inferior.

\vspace{-0.1in}
\subsection{Analysis of Combined Cycle Power Plant Data}

We also analyzed a real data set, namely the Combined cycle power plant data, for further evaluation of our proposed methods. This data set is available at the UCI Machine Learning Repository \url{https://archive.ics.uci.edu/ml/datasets/}. It contains 9568 observations from a Combined Cycle Power Plant over a period of six years (2006-2011), when the plant was set to work with full load. Each observation consists of hourly average values of ambient temperature, ambient pressure, relative humidity, exhaust vacuum and electric energy output. The idea was to predict electric energy output, which is dependent on other variables. When we used different methods to test for the independence among these five variables, all tests rejected $H_0$ on almost all occasions even when they were used on random subsets of size $10$ drawn from the data set. So, next we removed electric energy output from our analysis and carried out tests for independence among the other four variables.

\begin{figure}[h]
	\vspace{-0.1in}
	\centering
	\includegraphics[width=0.65\textwidth]{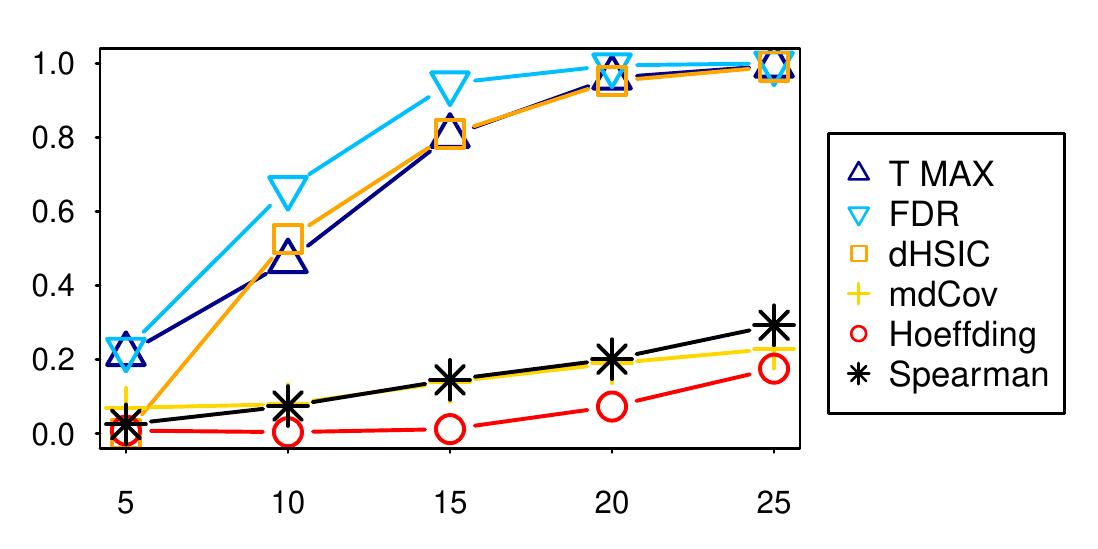}
	\vspace{-0.1in}
	\caption{Powers of different tests in  real datasets.}
	\label{fig: toird}
	\vspace{-0.2in}
\end{figure}

When we used the whole data set (ignoring electric energy output) for testing, all tests rejected $H_0$. It gives us a clear indication that these four variables have significant dependence among themselves, and different tests can be compared based on their powers. But based on a single experiment with the whole data set, it was not possible to compare among the powers of different test procedures. So, following the idea of \cite{sarkar2018some}, we carried out our experiments with subsets of different sizes. For each subset size (i.e., sample size), the experiment was repeated 10000 times to estimate the powers of different tests by proportion of times they rejected $H_0$. These estimated powers for different tests are shown in Figure \ref{fig: toird}. This figure clearly shows that in this example, our proposed test based on FDR outperformed its all competitors. The test based on $T_{Max}$ also performed well. Among the rest, only the dHSIC test had satisfactory performance.

\vspace{-0.1in}
\section{Discussion}

In this article, we have proposed a copula based multivariate dependency measure and established some of its theoretical properties. Unlike many other existing copula based measures, our dependency measure is invariant under strictly monotone transformations of the coordinate variables. Interestingly, it takes the value $0$ when the coordinate variables are independent and takes the value $1$ when for each pair of the coordinate variables, one is a strict monotone function of the other.  A data based estimate of this measure is proposed and some distribution-free tests of independence are constructed based on this estimate. Some nice theoretical properties of this estimate have also been derived and the large sample consistency of the resulting tests has been proved under appropriate regularity conditions. Unlike the dHSIC test, our proposed tests can be used even when the sample size is smaller than the number of variables. However, our proposed methods are not above all limitations. These rank based methods are mainly applicable when the coordinate variables are continuous in nature. In the case of discrete data, one may need to resolve the ties arbitrarily to define the ranks. The choice of the bandwidth is another issue to be resolved. In this article, we have adopted a multi-scale approach, where the results for different bandwidths are aggregated judiciously. The resulting tests worked well in all simulated and real data sets considered in this article. But instead of taking such a multi-scale approach, if we can choose a suitable data driven estimate of the bandwidth, that can further improve the performance of our methods, both in terms of power and computing time.
This can be considered as a problem for future research.

\begin{acknowledgements}
We would like to thank the Associate Editor and two anonymous reviewers for carefully reading the earlier version of the manuscript and also for providing us with several helpful comments. We are also thankful to Dr. Pierre Lafaye de Micheaux for providing the codes of the mdCov test. 
\end{acknowledgements}

\section*{Appendix}

\begin{proof}[{\bf Proof of Lemma 1}]
	For any permutation $\tau$ on $\mathbb{R}^d$, we have $\ks(\tau(\xvec),\tau(\yvec))$ $=\ks(\xvec,\yvec)$ and also, $\Tvec\sim\Pi$ implies $\tau(\Tvec)\sim\Pi$. Using these, one gets
	\vspace{-5pt}
	\begin{align*}
	\E_{(\Svec,\Svec^{'})\sim \mathsf{C}_{\tau(\Xvec)}\otimes\mathsf{C}_{\tau(\Xvec)}}[k_{\sigma}(\Svec,\Svec^{'})]=\E_{(\Svec,\Svec^{'})\sim \mathsf{C}_{\Xvec}\otimes\mathsf{C}_{\Xvec}}[k_{\sigma}(\tau(\Svec),\tau(\Svec^{'}))]\\
	=\E_{(\Svec,\Svec^{'})\sim \mathsf{C}_{\Xvec}\otimes\mathsf{C}_{\Xvec}}[k_{\sigma}(\Svec,\Svec^{'})]\; \mbox{ and}
	\end{align*}
		\vspace{-20pt}
	\begin{align*}
	\E_{(\Svec,\Tvec)\sim \mathsf{C}_{\tau(\Xvec)}\otimes \Pi}[k_{\sigma}(\Svec,\Tvec)]=\E_{(\Svec,\Tvec)\sim \mathsf{C}_{\Xvec}\otimes\Pi}[k_{\sigma}(\tau(\Svec),\tau(\Tvec))]\\
	=\E_{(\Svec,\Tvec)\sim \mathsf{C}_{\Xvec}\otimes\Pi}[k_{\sigma}(\Svec,\Tvec)].
	\end{align*}
	\vspace{-5pt}
	It follows that $\gamma_{\ks}(\mathsf{C}_{\tau(\Xvec)}, \Pi)=\gamma_{\ks}(\cx,\Pi)$ and hence $I_\sigma(\tau(\Xvec))=\cb$.\\
	
	Next, let $g:\mathbb{R}^d\rightarrow\mathbb{R}^d$ be a function of the form $g(x_1,x_2\cdots,x_d)=(g_1 (x_1)$, $g_2 (x_2),\cdots,g_d(x_d))$, where $g_{i_1},\cdots,g_{i_s}$ are strictly increasing and $g_{j_1},\cdots,g_{j_t}$ are strictly decreasing with $s+t=d$. Consider the function $f:\mathbb{R}^d\rightarrow\mathbb{R}^d$ given by $f(x_1,x_2,\cdots,x_d)=(f_1(x_1),f_2(x_2),\cdots,f_d(x_d))$, where $f_{i_l}(x)=x \;\forall\; l=1,2,\cdots,s$ and $f_{j_l}(x)=1-x \;\forall\; l=1,2,\cdots,t$. It can be easily verified that if $\Svec\sim \mathsf{C_{g(\Xvec)}}$ then $f(\Svec) \sim \mathsf{C_{\Xvec}}$. Applying this and the fact that $\ks(\Svec,\Svec^{'})=\ks(f(\Svec),f(\Svec^{'}))$, we get
	\vspace{-5pt}
	\begin{align*} \E_{(\Svec,\Svec^{'})\sim\mathsf{C}_{g(\Xvec)}\otimes\mathsf{C}_{g(\Xvec)}}[k_{\sigma}(\Svec,\Svec^{'})]=\E_{(\Svec,\Svec^{'})\sim\mathsf{C}_{g(\Xvec)}\otimes\mathsf{C}_{g(\Xvec})}[k_{\sigma}(f(\Svec),f(\Svec^{'}))]\\
	=\E_{(\Svec,\Svec^{'})\sim\mathsf{C}_{\Xvec}\otimes\mathsf{C}_{\Xvec}}[k_{\sigma}(\Svec,\Svec^{'})].
	\end{align*}
	By similar argument and using the fact that $\Tvec\sim\Pi$ implies $f(\Tvec)\sim\Pi$, one gets
	\begin{align*}\E_{(\Svec,\Tvec)\sim\mathsf{C}_{g(\Xvec)\otimes\Pi}}[k_{\sigma}(\Svec,\Tvec)]=\E_{(\Svec,\Tvec)\sim \mathsf{C}_{g(\Xvec)}\otimes\Pi}[k_{\sigma}(f(\Svec),f(\Tvec))]\\
	=\E_{(\Svec,\Tvec)\sim \cx\otimes\Pi}[k_{\sigma}(\Svec,\Tvec)].
	\end{align*}
	Thus $\gamma_{\ks}(\mathsf{C}_{g(\Xvec)}, \Pi)=\gamma_{\ks}(\mathsf{C}_{\Xvec},\Pi)$, whence, $I_\sigma(g(\Xvec))=\cb$, proving the invariance of $\cb$ under strictly monotonic transformations of $X_1,X_2,\ldots,X_d$. \qed
\end{proof}

\begin{proof}[{\bf Proof of Lemma 2}]
	Let $\Xvec$ be a random vector with continuous marginals, for which there is a $j$ such that each $X_i,\, i\not= j$ is a strictly monotonic function of $X_j$. Then, by Lemma 1, we have $\cb=I_\sigma(\Yvec)$ where $\Yvec=(X_j,X_j,\cdots,X_j)$. But then $\mathsf{C}_{\Yvec}$ is the maximum copula $\mathsf{M},\,$ so that, by definition, $I_\sigma(\Yvec)=1$. \qed
\end{proof}

\begin{proof}[{\bf Proof of Theorem 1}]
	This proof has two steps. At the first step, we prove that $\csd=\gamma^2_{\ks}(\mathsf{M},\Pi)$. At the second step, we prove that $\gamma^2_{\ks}(\cx,\Pi)=\left(\frac{\sigma}{\sqrt{2\pi}}\right)^d \int_{\mathbb{R}^d}|\varphi_{\cx}(\wvec)-\varphi_{\Pi}(\wvec)|^2\exp\left(-\frac{\sigma^2}{2}\wvec^{\intercal}\wvec\right)\,d\wvec$. Clearly, proving these two steps will complete the proof.
	
	\noindent
	\underline{First  step:}
	Note that for $(\Svec,\Svec^{'}\!\!, \Tvec, \Tvec^{'})\sim \mathsf{M}\otimes\mathsf{M}\otimes\Pi\otimes\Pi$, we have
	\vspace{-5pt}
	\begin{align*}
	&\gamma^2_{\ks}(\mathsf{M},\Pi)=\E[k_{\sigma}(\Svec,\Svec^{'})]-2\E[\ks(\Svec,\Tvec)]+\E[k_{\sigma}(\Tvec,\Tvec^{'})]\\
	&=\int^1_0\!\!\!\int^1_0\!\! e^{-\frac{d(u-v)^2}{2\sigma^2}}\,du\,dv-2\int^1_0\!\!\left[\int^1_0\!\! e^{-\frac{(u-v)^2}{2\sigma^2}}\,du\right]^d \,dv+\left[\int^1_0\!\!\!\int^1_0 e^{-\frac{(u-v)^2}{2\sigma^2}}\,du\,dv\right]^d\\
	&=\kappa\left(\frac{\sigma}{\sqrt{d}}\right)-2\int^1_0\lambda^d(u,\sigma)\,du+\kappa^d(\sigma)=\csd.
	\end{align*}
	
	\vspace{0.025in}
	\noindent
	\underline{Second step:} 
	We use the well-known formula for Fourier transform of the\\ $d$-dimensional Gaussian density:
	\vspace{-5pt}
	$$\exp\left(-\frac{1}{2\sigma^2}\xvec^{\intercal}\xvec\right)=\int_{\mathbb{R}^d}e^{-\sqrt{-1}\xvec^{\intercal}\wvec}\cdot\left(\frac{\sigma}{\sqrt{2\pi}}\right)^d\exp\left(-\frac{\sigma^2}{2}\wvec^{\intercal}\wvec\right)\,d\wvec,\; \xvec\in\mathbb{R}^d.$$
	This gives us
	\vspace{-5pt}
	$$ \ks(\xvec,\yvec)=\left(\frac{\sigma}{\sqrt{2\pi}}\right)^d \int_{\mathbb{R}^d}e^{-\sqrt{-1}\xvec^{\intercal}\wvec}\cdot e^{\sqrt{-1}\yvec^{\intercal}\wvec}\exp\left(-\frac{\sigma^2}{2}\wvec^{\intercal}\wvec\right)\,d\wvec,\; \xvec,\yvec\in\mathbb{R}^d.$$
	Using the representation of $\gamma^2_{k}$ from equation (2) and Fubini's theorem, one gets
	\vspace{-5pt}
	\begin{align*}
	\gamma^2_{\ks}(\cx,\Pi)&=\!\left(\frac{\sigma}{\sqrt{2\pi}}\right)^{\!d}\!\!\int_{\mathbb{R}^d}\Big[\,\varphi_{\cx}(\wvec)\overline{\varphi_{\cx}}(\wvec)+\varphi_{\Pi}(\wvec)\overline{\varphi_{\Pi}}(\wvec)\\
	&\qquad\qquad\qquad-2\varphi_{\cx}(\wvec)\overline{\varphi_{\Pi}}(\wvec)\,\Big]\,\exp\left(-\frac{\sigma^2}{2}\wvec^{\intercal}\wvec\right)\,d\wvec,
	\end{align*}
	from which the second part follows. \qed
\end{proof}

\begin{customlemma}{L1}
	\label{lem: cov form}
	Let $(\Xvec,\Yvec)$ and $(\Xvec^{'},\Yvec^{'})$ be independent and identically distributed random vectors taking values in $\mathcal{\Xvec}\times\mathcal{\Yvec}$. Given symmetric measurable functions $k:\mathcal{X}\times\mathcal{X}\to\mathbb{R}$ and $\overline{k}:\mathcal{Y}\times\mathcal{Y}\to\mathbb{R}$, define \vspace{-1ex}
	\begin{align*}
	&V=k(\Xvec,\Xvec^{'})-\E\left[k(\Xvec,\Xvec^{'})\Big|\Xvec\right]-\E\left[k(\Xvec,\Xvec^{'})\Big|\Xvec^{'}\right]+\E\left[ k(\Xvec,\Xvec^{'})\right]\\
	&W=\overline{k}(\Yvec,\Yvec^{'})-\E\left[\,\overline{k}(\Yvec,\Yvec^{'})\Big|\Yvec\right]-\E\left[\,\overline{k}(\Yvec,\Yvec^{'})\Big|\Yvec^{'}\right]+\E\left[\, \overline{k}(\Yvec,\Yvec^{'})\right].
	\end{align*}\vspace{-1ex}
	Then, we have \vspace{-1ex}
	\begin{align*}
	\E\left[\,VW\,\right]=\,\E\left[k(\Xvec,\Xvec^{'})\,\overline{k}(\Yvec,\Yvec^{'})\right]+\E\left[k(\Xvec,\Xvec^{'})\right]\E\left[\overline{k}(\Yvec,\Yvec^{'})\right]\\
	-2\,\E\left[\, \E\left[k(\Xvec,\Xvec^{'})\Big|\Xvec\right]\E\left[\overline{k}(\Yvec,\Yvec^{'})\Big|\Yvec\right]\,\right].
	\end{align*}
\end{customlemma}

\begin{proof}
	The proof is based on expanding the product $VW$ and then taking term-by-term expectations. One and only one term gives $\E\left[k(\Xvec,\Xvec^{'})\,\overline{k}(\Yvec,\Yvec^{'})\right]$. The seven terms, where at least one of $\E\left[k(\Xvec,\Xvec^{'})\right]$ or
	$\E\left[\overline{k}(\Yvec,\Yvec^{'})\right]$ appear as a factor, and the two terms $\E\left[k(\Xvec,\Xvec^{'})\Big|\Xvec\right]\cdot \E\left[\,\overline{k}(\Yvec,\Yvec^{'})\Big|\Yvec^{'}\right]$ and $\E\left[k(\Xvec,\Xvec^{'})\Big|\Xvec^{'}\right]\cdot \E\left[\,\overline{k}(\Yvec,\Yvec^{'})\Big|\Yvec\right]$, will all give the same expectation $\E\left[k(\Xvec,\Xvec^{'})\right]\E\left[\overline{k}(\Yvec,\Yvec^{'})\right]$ (the last two because of independence of $(\Xvec,\Yvec)$ and $(\Xvec^{'},\Yvec^{'})$). Taking into account the signs of these nine terms with the same expectation, we would be left with just one with a positive sign. Next, the remaining six terms will all have the same expectation, namely, $\E\left[ \E\left[k(\Xvec,\Xvec^{'})\Big|\Xvec\right]\E\left[\overline{k}(\Yvec,\Yvec^{'})\Big|\Yvec\right]\right]$.  For two of the terms, this is straightforward. But the other four terms need judicious use of properties of conditional expectation. For example, by independence of $(\Xvec,\Yvec)$ and $(\Xvec^{'},\Yvec^{'})$, we have $\E\left[\overline{k}(\Yvec,\Yvec^{'})\Big|\Yvec\right]=\E\left[\overline{k}(\Yvec,\Yvec^{'})\Big|(\Xvec,\Yvec)\right]$ and similarly $\E\left[k(\Xvec,\Xvec^{'})\Big|(\Xvec,\Yvec)\right]=\E\left[k(\Xvec,\Xvec^{'})\Big|\Xvec\right]$. Using these, we get \vspace{-1.5ex}
	\begin{align*}
	\E\left[\,k(\Xvec,\Xvec^{'})\E\left[\overline{k}(\Yvec,\Yvec^{'})\Big|\Yvec\right]\,\right]&=
	\E\left[\,k(\Xvec,\Xvec^{'})\E\left[\overline{k}(\Yvec,\Yvec^{'})\Big|(\Xvec,\Yvec)\right]\,\right]\\
	&=
	\E\left[\,\E\left[k(\Xvec,\Xvec^{'})\big|(\Xvec,\Yvec)\right]\E\left[\overline{k}(\Yvec,\Yvec^{'})\Big|(\Xvec,\Yvec)\right]\,\right]\\
	&=
	\E\left[\,\E\left[k(\Xvec,\Xvec^{'})\big| \Xvec\right]\E\left[\overline{k}(\Yvec,\Yvec^{'})\Big| \Yvec\right]\,\right].
	\end{align*}
	One can similarly handle other three terms. Considering the signs of these six terms with the same expectation, one is left with $-2\E\left[\E\left[k(\Xvec,\Xvec^{'})\big|\Xvec\right]\E\left[\overline{k}(\Yvec,\Yvec^{'})\Big| \Yvec\right]\right]$. This completes the proof. \qed
\end{proof}

\begin{proof}[{\bf Proof of Theorem 2}]
	(a) By definition, $V_1$ and $V_2$ have zero means. Using Lemma \ref{lem: cov form} with kernels $k=\overline{k}=\ks$ on
	$\mathcal{X}=\mathcal{Y}=\mathbb{R}^d$, we get
	\vspace{-5pt}
	\begin{align*}
	\cov[V_1,V_2]=&\E[V_1V_2]=\E\left[k_{\sigma}(T_1,T_1^{'})\ks(T_2,T_2^{'})\right]+\E\left[\ks(T_1,T_1^{'})\right]\E\left[\ks(T_2,T_2^{'})\right]\\
	&-2\,\E\left[ \E\left[k_{\sigma}(T_1,T_1^{'})\Big|T_1\right]\E\left[\ks(T_2,T_2^{'})\Big|T_2\right]\right]=\gamma^2_{\ks}(\mathsf{C}_{(X,Y)},\Pi).
	\end{align*}
	One can similarly show that $\var[V_1]=\var[V_2]=\gamma^2_{\ks}(\mathsf{M},\Pi)$ and hence $I^2_\sigma(X,Y)=\cor[V_1,V_2]$.
	The inequality $I_\sigma(X_1,X_2)\leq 1$  follows from it. Further, from the condition for equality in Cauchy-Schwartz inequality and the fact that $V_1$ and $V_2$ are identically distributed, it follows that $I_\sigma(X_1,X_2)=1$ if and only if $V_1=V_2$ almost surely.
	
	Since $T_1$ and $T_1^{'}$ are independent and uniformly distributed random variables on $[0,1]$ and so also are $T_2$ and $T_2^{'}$, it follows that $V_1=V_2$ almost surely if and only if $g(T_1,T_1^{'})=g(T_2,T_2^{'})$ almost surely, where $g(x,y)=\ks(x,y)-\lambda(x,\sigma )-\lambda(y,\sigma )$ with $\lambda (\cdot ,\sigma )$ as defined in Theorem 1.
	
	Now, using the facts that $(T_1,T_2)$ and $(T_1^{'},T_2^{'})$ are independent and identically
	distributed with values in $[0,1]^2$ and that the function $g$ is uniformly continuous on the compact set $[0,1]^2$, one can easily deduce that $g(T_1,T_1^{'})=g(T_2,T_2^{'})~\text{a.s.}$ implies $g(T_1,T_1)=g(T_2,T_2)$~\text{a.s.} But, this, in turn, implies that $\Phi\left(\frac{T_1}{\sigma}\right)+\Phi\left(\frac{1-T_1}{\sigma}\right)=\Phi\left(\frac{T_2}{\sigma}\right)+\Phi\left(\frac{1-T_2}{\sigma}\right)~\text{a.s.}$. From this, we may conclude that $\pr\left[ T_2=T_1~\text{or } T_2=1-T_1\right]=1$ and also $\pr\left[\lambda(T_1,\sigma )=\lambda(T_2,\sigma )\right]=1$. Of course, the same would be true of the pair $(T_1^{'},T_2^{'})$, which is moreover independent of the pair $(T_1,T_2)$.
	
	Using these in the equality $g(T_1,T_1^{'})=g(T_2,T_2^{'})~\text{a.s.}$, one obtains $\ks(T_1,T_1^{'})=\ks(T_2,T_2^{'})~\text{a.s.}$, which implies that
	$|T_1-T_1^{'}|=|T_2-T_2^{'}|~\text{a.s.}$ We conclude that either $T_2=T_1~\text{a.s.}$ or $T_2=1-T_1~\text{a.s.}$ Thus the copula distribution of $(X_1,X_2)$ is
	either the distribution of $(T_1,T_1)$ or that of $(T_1,1-T_1)$ where $T_1$ is uniformly distributed on $[0,1]$. So, 
	$X_1$ and $X_2$ are almost surely strictly monotone functions of each other.
	
	\vspace{0.05in}
	\noindent(b)
	For $|r|<1$, let $\phi_{r}$ denote the density of the standard bivariate normal distribution with correlation coefficient $r$. Also, let $\Phi$ and $\phi$ denote respectively the cumulative distribution function and the density function of the standard univariate normal distribution. It is well-known that the copula distribution of any bivariate normal distribution with correlation coefficient $r$ is the same as that of the standard bivariate normal distribution with the same correlation coefficient. Using the well-known Mehler's representation (see  \cite{kibble1945}, Page 1) of standard bivariate normal density with correlation $r$, one then gets that, for $|r|<1$, the copula distribution $\mathcal{C}(r)$ of any bivariate normal distribution with correlation coefficient $r$ has density given by
	\vspace{-5pt}
	$$ \frac{\phi_{r}(\Phi^{-1}(u),\Phi^{-1}(v))}{\phi(\Phi^{-1}(u))\phi(\Phi^{-1}(v))}=\sum_{i=0}^{\infty}\frac{r^i}{i!}H_i((\Phi^{-1}(u))H_i(\Phi^{-1}(v)), \quad (u,v)\in [0,1]^2,\vspace{-5pt} 
	$$
	where $\{H_i(x), i\geq 0\}$ are the well-known Hermite polynomials.
	Using this, we get that if $(S,T)\sim \mathcal{C}(r_1)\otimes \mathcal{C}(r_2)$ with $|r_1|<1,|r_2|<1$, then
	\vspace{-5pt}
	\begin{align}
	\E[\ks(S,T)]
	\!=&\!\!\!\int\limits_{[0,1]^4}\! e^{-\frac{(s_1-t_1)^2+(s_2-t_2)^2}{2\sigma^2}}\sum_{i=0}^{\infty}\frac{r_1^i}{i!} H_i(\Phi^{-1}(s_1))H_i(\Phi^{-1}(s_2))\nonumber\\
	&\times \sum_{j=0}^{\infty}\!\frac{r_2^j}{j!}H_j(\Phi^{-1}(t_1))H_j(\Phi^{-1}(t_2))\,ds_1 \,ds_2 \,dt_1 \,dt_2\vspace{-10pt}\tag{E1}\label{eq: E1}
	\end{align}
	
	We now claim that in the above expression, the double summation and integration can be interchanged. To justify this, we recall that the Hermite polynomials $\{H_i(\cdot ), i\geq 0\}$ form a complete orthonormal basis for $L_2(\mathbb{R},\phi(x)dx)$ and, in particular,  for any $i\geq 0$,  $\int_{[0,1]}\left|H_i(\Phi^{-1}(s))\right|\, ds=
	\int_{\mathbb{R}}\left| H_i(x)\right| \phi(x)\,dx\leq \left[\int_{\mathbb{R}}
	H_i^2(x)\phi(x)\, dx\right]^{\frac{1}{2}}$ $=1$. As a consequence,  \vspace{-.5ex}
	\begin{align*}
	&\sum_{i=0}^{\infty}\sum_{j=0}^{\infty}\int_{[0,1]^4}\left|e^{-\frac{(s_1-t_1)^2+(s_2-t_2)^2}{2\sigma^2}}\right|\,
	\left|\frac{r_1^ir_2^j}{i!j!}\right|\,\left| H_i(\Phi^{-1}(s_1))\right|\,\left|H_i(\Phi^{-1}(s_2))\right|\\
	&\qquad\times \left|H_j(\Phi^{-1}(t_1))\right|\left|H_j(\Phi^{-1}(t_2))\right|\,ds_1 \,ds_2\, dt_1 \,dt_2\\
	&\leq\sum_{i=0}^{\infty}\sum_{j=0}^{\infty}\frac{1}{i!j!}\left[\int\limits_{[0,1]}\!\!\left| H_i(\Phi^{-1}(s))\right|ds\right]^2\!\left[\int\limits_{[0,1]}
	\!\!\left|H_j(\Phi^{-1}(t))\right|dt\right]^2 \!\!\leq\!\sum_{i=0}^{\infty}\sum_{j=0}^{\infty}\frac{1}{i!j!}<\infty .
	\end{align*}
	
	We can, therefore, interchange the double summation and integration in the right-hand-side of the equation \eqref{eq: E1} above to obtain that, for any $r_1,r_2$ with $|r_1|<1,|r_2|<1$,
		\vspace{-5pt}\begin{align*}
	&\E_{(S,T)\sim \mathcal{C}(r_1)\otimes \mathcal{C}(r_2)}[\ks(S,T)]=
	\sum_{i=0}^{\infty}\sum_{j=0}^{\infty}a_{i,j}r_1^ir_2^j\, ,
	\end{align*}
	\vspace{-0.2in}
	\begin{align}
	\label{eq: E2}
	{\mbox where}~~ a_{i,j}&:=\frac{1}{i!j!}\left[\int_{[0,1]^2} e^{-\frac{(u-v)^2}{2\sigma^2}}H_i(\Phi^{-1}(s))H_j(\Phi^{-1}(t))\,ds \,dt\right]^{\!2}\nonumber\\
	&=\frac{1}{i!j!}\left[\int_{\mathbb{R}^2} e^{-\frac{(\Phi(x)-\Phi(y))^2}{2\sigma^2}}H_i(x)H_j(y)\phi(x)\phi(y)\,dx \,dy\right]^{\!2}.\tag{E2}
	\end{align}
	Observe that $a_{i,j}\geq 0, a_{i,j}=a_{j,i}$ and also, $(i!j!)a_{i,j}\!\leq\!\left[\int\limits_{\mathbb{R}}\left| H_i(x)\right| \phi(x)\,dx\right]^2\!\leq\! 1$. \\
	Note that for any bivariate normal random vector $(X_1,X_2)$ with correlation coefficient $r$ (where $|r|<1$), we have
	$\gamma^2_{\ks}(\mathcal{C}_{(X_1,X_2)}, \Pi)=\gamma^2_{\ks}(\mathcal{C}(r),\mathcal{C}(0))$, which equals
	\begin{align*}
	&\E_{_{(S,S^{'})\sim \mathcal{C}(r)\otimes\mathcal{C}(r)}}\![k_{\sigma}(S,\!S^{'})]\!
	-\!2\E_{_{(S,T)\sim \mathcal{C}
			(r)\otimes \mathcal{C}(0)}}\![k_{\sigma}(S,\!T)]
	\!+\!\E_{_{(T,T^{'})\sim\mathcal{C}(0)\otimes\mathcal{C}(0)}}\![k_{\sigma}(T,\!T^{'})]\\
	&=\!\sum_{k=0}^{\infty}\!\sum_{\substack{i,j\geq 0 \\ i+j=k}}\!a_{i,j}\!r^k \!-\!2\,\sum_{k=0}^{\infty}a_{k,0}\,r^k\!+\!a_{0,0}
	=\!\sum_{k=1}^{\infty}\!\sum_{\substack{i,j \geq 1\\ i+j=k}}\!a_{i,j}\;r^k\underset{(a)}{=}\,
	\sum_{k=1}^{\infty}\!\sum_{\substack{i,j\geq 1\\ i+j=2k}}\!a_{i,j}\;r^{2k}.
	\end{align*}
	Equality (a) is due to the fact that the $i$th Hermite polynomial $H_i$ is an even or an odd function according as $i$ is even or odd, so that if exactly one of $i$ and $j$ is odd, then $a_{i,j}=0$, as can easily be seen from equation \eqref{eq: E2}. Therefore, we have
	\vspace{-1ex}
	\begin{align*}
	\label{eq: rho}
	I^2_\sigma(X_1,X_2)&=\gamma^{-2}_{\ks}(\mathsf{M},\Pi)\gamma^2_{\ks}(\mathcal{C}_{(X_1,X_2)},\Pi)=\gamma^{-2}_{\ks}(\mathsf{M},\Pi)\sum_{k=1}^{\infty}\sum_{\substack{i,j\geq 1, i+j=2k}}\!\!a_{i,j}\;r^{2k}.
	\vspace{-10pt}
	\end{align*}
	So, 	$I^2_\sigma(X_1,X_2)=r^2g(r)$, where $g(r)=\gamma^{-2}_{\ks}(\mathsf{M},\Pi)\sum\limits_{k=1}^{\infty}\;\sum\limits_{i,j\geq 1:\, i+j=2k}\!\!a_{i,j}\;r^{2(k-1)}$ is a power series in $r^2$ with positive coefficients and hence increasing in $|r|$. So, $I^2_\sigma(X_1,X_2)=r^2\cdot g(r)$ is an increasing function of $|r|$. \qed
\end{proof}

\begin{proof}[{\bf Proof of Theorem 3}]
	It is enough to show that for every dimension $d~ (\geq 3)$, there exist two $d$ dimensional copulas $C_1$ and $C_2$ with $\mathcal{M}(C_1)\neq\mathcal{M}(C_2)$, such that for any choice of co-ordinates $\{ i_1,\cdots , i_k\}\subsetneqq\{1,\ldots ,d\}$, if $C_1^{'}$ and $C_2^{'}$ are the associated marginal copulas arising out of of $C_1$ and $C_2$, then $\mathcal{M}(C_1^{'})=\mathcal{M}(C_2^{'})$.
	
	Take $C_1$ to be the $d$-dimensional uniform copula $\Pi$. Then $\mathcal{M}(C_1)=0$, and also for any lower dimensional marginal copula $C_1^{'}$ of $C_1$, $\mathcal{M}(C_1^{'})=0$. We now exhibit a $d$-dimensional copula $C_2\neq\Pi$ such that any lower dimensional marginal copula $C_2^{'}$ of $C_2$ is uniform copula. We would then have $\mathcal{M}(C_1)=0\neq\mathcal{M}(C_2)$ but $\mathcal{M}(C_1^{'})=\mathcal{M}(C_2^{'})=0$, which will complete the proof.
	
	We take $C_2$ to be the copula given by the copula density $\mathcal{C}_2$ defined as
	$$\mathcal{C}_2(u_1,u_2,,\cdots,u_d)=2\,\mathbb{I}\left[ \left(u_1-\frac{1}{2}\right) \left(u_2-\frac{1}{2}\right) \cdots \left(u_d-\frac{1}{2}\right) \geq 0\right],$$
	where $\mathbb{I}$ denotes the indicator function. To show that all lower dimensional marginal copulas of $C_2$ are uniform, it is enough to show that the marginal copula $C_2^{'}$ that we get from $C_2$ discarding the $d^{\text{th}}$ co-ordinate, is uniform. Now, note that the density of $C_2^{'}$ is given by
	\vspace{-10pt}
	\begin{align*}
	\mathcal{C}^{'}_2(u_1,u_2,,\cdots,u_{d-1})&=\int_{0}^{1}2\mathbb{I}\left[ \left(u_1-\frac{1}{2}\right) \left(u_2-\frac{1}{2}\right) \cdots \left(u_d-\frac{1}{2}\right) \geq 0\right]\,du_d\\
	&=\int_{0}^{\frac{1}{2}}2\mathbb{I}\left[ \left(u_1-\frac{1}{2}\right) \left(u_2-\frac{1}{2}\right) \cdots \left(u_{d-1}-\frac{1}{2}\right) \leq 0\right]\,du_d\\
	&\quad+\int_{\frac{1}{2}}^{1}2\mathbb{I}\left[ \left(u_1-\frac{1}{2}\right) \left(u_2-\frac{1}{2}\right) \cdots \left(u_{d-1}-\frac{1}{2}\right) \geq 0\right]\,du_d\\
	&=\mathbb{I}\left[ \left(u_1-\frac{1}{2}\right) \left(u_2-\frac{1}{2}\right) \cdots \left(u_{d-1}-\frac{1}{2}\right) \leq 0\right]\\
	&\quad+\mathbb{I}\left[ \left(u_1-\frac{1}{2}\right) \left(u_2-\frac{1}{2}\right) \cdots \left(u_{d-1}-\frac{1}{2}\right) \geq 0\right]=1. ~~~~\hfill \qed
	\end{align*}
\end{proof}

\begin{proof}[{\bf Proof of Lemma 4}]
	We shall prove that as $\sigma\to\infty$, $\sigma^4\gamma^2_{k_\sigma}(\cx,\Pi)\to\sum_{1\leq i< j\leq d}\cov^2(S_i,S_j)$. It will imply that  $\sigma^4\gamma^2_{k_\sigma}(\mathsf{M},\Pi)\to{d \choose 2}\var^2(S_1)$ as $\sigma\to\infty$, which  in turn will imply that  $\scb=\frac{\sigma^4\gamma^2_{k_\sigma}(\cx,\Pi)}{\sigma^4\gamma^2_{k_\sigma}(\mathsf{M},\Pi)}\to\frac{1}{{d \choose 2}}\sum_{1\leq i< j\leq d}\frac{\cov^2(S_i,S_j)}{\var^2(S_1)}=\frac{1}{{d \choose 2}}\sum_{1\leq i< j\leq d}\frac{\cov^2(S_i,S_j)}{\var(S_i)\var(S_j)}=\frac{1}{{d \choose 2}}\sum_{1\leq i< j\leq d}\cor^2(S_i,S_j)$ as $\sigma\to\infty$; which is our desired result. 
	
	Observe that $\E k_\sigma(\Svec,\Tvec)=1-\frac{1}{2\sigma^2}\E\|\Svec-\Tvec\|^2_2+\frac{1}{8\sigma^4}\E\|\Svec-\Tvec\|^4_2+\mathcal{O}\left(\frac{1}{\sigma^6}\right).$ Assume that $\Svec,\Svec^{'},\Tvec$ and $\Tvec^{'}$ are four random vectors such that $(\Svec,\Svec^{'},\Tvec,\Tvec^{'})\sim\cx\otimes\cx\otimes\Pi\otimes\Pi$. Then,
	\begin{align*}
	\gamma^2_{k_\sigma}(\cx,\Pi)&=\E k_\sigma(\Svec,\Svec^{'})-2\E k_\sigma(\Svec,\Tvec)+\E k_\sigma(\Tvec,\Tvec^{'})\end{align*} \begin{align*}
	&=-\frac{1}{2\sigma^2}\E\left[\|\Svec-\Svec^{'}\|_2^2+\|\Tvec-\Tvec^{'}\|_2^2-2\|\Svec-\Tvec\|_2^2\right]\\
	&\quad+\frac{1}{8\sigma^4}\E\left[\|\Svec-\Svec^{'}\|_2^4+\|\Tvec-\Tvec^{'}\|_2^4-2\|\Svec-\Tvec\|_2^4\right]+\mathcal{O}\left(\frac{1}{\sigma^6}\right).
	\end{align*}
	\vspace{-10pt}
	\begin{align*}
	\mbox{Now, }&\E\left[\|\Svec-\Svec^{'}\|_2^2+\|\Tvec-\Tvec^{'}\|_2^2-2\|\Svec-\Tvec\|_2^2\right]\\
	&\quad\quad\quad\quad\quad\quad\quad=\E\sum_{i=1}^{d}\left[(S_i-S_i^{'})^2+(T_i-T_i^{'})^2-2(S_i-T_i)^2\right]=0, \mbox{ and }
	\end{align*} \vspace{-15pt}
	\begin{align*}
	&\E\left[\|\Svec-\Svec^{'}\|_2^4+\|\Tvec-\Tvec^{'}\|_2^4-2\|\Svec-\Tvec\|_2^4\right]\\
	&=\E\sum_{i=1}^{d}\left[(S_i-S_i^{'})^4+(T_i-T_i^{'})^4-2(S_i-T_i)^4\right]\\
	&+2\E\hspace{-0.15in}\sum_{1\leq i< j\leq d}\hspace{-0.15in}\Big[(S_i-S_i^{'})^2(S_j-S_j^{'})^2+(T_i-T_i^{'})^2(T_j-T_j^{'})^2
	-2(S_i-T_i)^2(S_j-T_j)^2\Big]\\
	&=2\E\hspace{-0.15in}\sum_{1\leq i< j\leq d} \hspace{-0.15in}\Big[(S_i-S_i^{'})^2(S_j-S_j^{'})^2+(T_i-T_i^{'})^2(T_j-T_j^{'})^2
	-2(S_i-T_i)^2(S_j-T_j)^2\Big].
	\end{align*}
	Hence $\gamma^2_{k_\sigma}(\cx,\Pi)=\frac{1}{4\sigma^4}\E\sum_{1\leq i< j\leq d}\Big[(S_i-S_i^{'})^2(S_j-S_j^{'})^2+(T_i-T_i^{'})^2(T_j-T_j^{'})^2-2(S_i-T_i)^2(S_j-T_j)^2\Big]+\mathcal{O}\left(\frac{1}{\sigma^6}\right).$ Now, by some straight-forward  but tedious calculations, it can be shown that $\E\Big[(S_i-S_i^{'})^2(S_j-S_j^{'})^2+(T_i-T_i^{'})^2(T_j-T_j^{'})^2-2(S_i-T_i)^2(S_j-T_j)^2\Big]=4\cov^2(S_i,S_j)$. This implies that as $\sigma\to\infty$, $\sigma^4\gamma^2_{k_\sigma}(\cx,\Pi)\to\sum_{1\leq i< j\leq d}\cov^2(S_i,S_j).$ \qed
\end{proof}

\begin{proof}[{\bf Proof of Lemma 5}]
	(a)
	Clearly, applying a permutation to the coordinates of the observation vectors $\Xvec^{(i)}, i=1,2,\cdots , n,$ changes the coordinates of the $\Yvec^{(i)}$'s by the same permutation. Since $s_1$ and $s_2$ from equation (7) are both invariant under permutation of coordinates of the $\Yvec^{(i)}$'s, the proof is complete.
	
	\vspace{0.05in}
	\noindent
	For proving invariance under monotonic transformation, it is enough to consider the case when only one of the coordinates in the observation vectors is changed by a strictly monotonic non-identity transformation. Assume, therefore, that only the $s^{\text{th}}$ coordinate of the $\Xvec^{(i)}$'s is changed by a strictly monotonic transformation, while the other coordinates are kept the same. This will affect only the $s^{\text{th}}$ coordinate of the $\Yvec^{(i)}$'s. Denoting the changed $\Yvec^{(i)}$'s as $\Yvec^{*(i)}$'s, it is clear that $Y^{*(i)}_s$ will equal $Y^{(i)}_s$ or $1+\frac{1}{n}-Y^{(i)}_s$ for all $i=1,2,\cdots,n$, according as the transformation is strictly increasing or strictly decreasing. In either case, $\ks(\Yvec^{*(i)},\Yvec^{*(j)})=\ks(\Yvec^{(i)},\Yvec^{(j)})$, so that $s_1$ in equation (7) remains unchanged. One can easily see that $s_2$ also remains unchanged as well.
	
	\vspace{0.05in}
	(b)	Without loss of generality, we may assume that the first coordinates of the $\Xvec^{(i)}$'s are in ascending order. Now, suppose that every other coordinate of the $\Xvec^{(i)}$'s is in a strictly monotonic relation with the first coordinate; then, for $j=2,\cdots , d$, the $j^{\text{th}}$ coordinates of the $\Xvec^{(i)}$'s will be in either ascending or descending order. By monotonic transformation invariance property, we may assume, without loss of generality, that all the coordinates of the $\Xvec^{(i)}$'s are in ascending order. But then, the $\Yvec^{(i)}$'s are clearly given by $Y_j^{(i)}=\frac{i}{n}$, for all $j$ and one can then see that $s_1=v_1$ and $s_2=v_2$, whence it follows that $\ecb=1$. \qed
\end{proof}

\begin{proof}[{\bf Proof of Theorem 4}] For two independent random vectors
 $(T_1,T_2)$ and $(T_1^{'},T_2^{'})$  with both having distribution $\cxn$, one has \vspace{-1.5ex}	$$V_1=\ks(T_1,T_1^{'})-\E\left[\ks(T_1,T_1^{'})\Big|T_1\right]-\E\left[\ks(T_1,T_1^{'})\Big|T_1^{'}\right]+\E\left[\ks(T_1,T_1^{'})\right] \mbox { and}$$ \vspace{-4ex} $$V_2=\ks(T_2,T_2^{'})-\E\left[\ks(T_2,T_2^{'})\Big|T_2\right]-\E\left[\ks(T_2,T_2^{'})\Big|T_2^{'}\right]+\E\left[\ks(T_2,T_2^{'})\right] \mbox { so that}  \vspace{-1.5ex}$$
	\begin{align*}
	&\gamma^2_{\ks}(\cxn,\Pi_n)\\
	&=\E\left[k_{\sigma}(T_1,T_1^{'})\ks(T_2,T_2^{'})\right]-2\E\left[ \E\left[k_{\sigma}(T_1,T_1^{'})\Big|T_1\right]\E\left[\ks(T_2,T_2^{'})\Big|T_2\right]\right]\\
	&+\E\left[\ks(T_1,T_1^{'})\right]\E\left[\ks(T_2,T_2^{'})\right]\\
	&=\E\left[\left\{\ks(T_1,T_1^{'})-\E\left[\ks(T_1,T_1^{'})\Big|T_1\right]-\E\left[\ks(T_1,T_1^{'})\Big|T_1^{'}\right]+\E\left[\ks(T_1,T_1^{'})\right]\right\}\right.\\
	&\quad\quad\quad\left.\left\{\ks(T_2,T_2^{'})-\E\left[\ks(T_2,T_2^{'})\Big|T_2\right]-\E\left[\ks(T_2,T_2^{'})\Big|T_2^{'}\right]+\E\left[\ks(T_2,T_2^{'})\right]\right\}\right]\\
	&=\frac{1}{n^2}\sum_{1\leq i,j\leq n}V_1(i,j)V_2(i,j)\, .
	\end{align*}
	Here the second last equality follows from Lemma \ref{lem: cov form}. Similarly, one can show that
	\vspace{-1ex}
	$$\gamma^2_{\ks}(\mathsf{M}_n,\Pi_n)=\frac{1}{n^2}\sum_{1\leq i,j\leq n}V_1^2(i,j)=\frac{1}{n^2}\sum_{1\leq i,j\leq n}V_2^2(i,j)\, .$$
	
	Cauchy-Schwartz (CS) inequality immediately gives $\ecb^2\leq 1$. Further, by the necessary and sufficient condition for equality in the CS inequality and using the fact that $\sum_{1\leq i,j\leq n}V^2_1(i,j)=\sum_{1\leq i,j\leq n}V^2_2(i,j)$, one gets that $\ecb=1$ if and only if $V_1(i,j)=V_2(i,j)~\forall~i,j$.
	
	Now, if one coordinate of the observation vectors is a monotone function of the other coordinate, then either $Y_2^{(i)}=Y_1^{(i)}~\forall~ i$ or $Y_2^{(i)}=\frac{n+1}{n}-Y_1^{(i)}~\forall~ i$. In either case, $|Y_1^{(i)}-Y_1^{(j)}|=|Y_2^{(i)}-Y_2^{(j)}|~\forall~i, j$, which will clearly imply that $V_1(i,j)=V_2(i,j)~\forall~i,j$.
	
	To prove the converse, first observe that for any $i$,
	\begin{align*}
	&\sum_{l=1}^{n}\ks(Y_1^{(i)},Y_1^{(l)})=\sum_{l=1}^{n}\ks(Y_1^{(l)},Y_1^{(i)})=
	\sum_{l=1}^{n}\ks(Y_1^{(i)},l/n) \text{~and~}\\
	&\sum_{l=1}^{n}\ks(Y_2^{(i)},Y_2^{(l)})=\sum_{l=1}^{n}\ks(Y_2^{(l)},Y_2^{(i)})=\sum_{l=1}^{n}\ks(Y_2^{(i)},l/n).
	\end{align*}
	Now suppose that $V_1(i,j)=V_2(i,j)~\forall~i,j$. Then, taking $i=j$, one deduces that
	\begin{equation}
	\sum_{l=1}^{n} \ks(Y^{(i)}_1,l/n)=\sum_{l=1}^{n} \ks(Y^{(i)}_2,l/n)\quad\forall~i\in\{1,2,\cdots,n\}.\label{eq: E4}\tag{E4}
	\end{equation}
	Using this now in $V_1(i,j)=V_2(i,j)$, one gets
	\begin{equation}
	\ks(Y_1^{(i)},Y_1^{(j)})=\ks(Y_2^{(i)},Y_2^{(j)}),\text{ i.e., } |Y_1^{(i)}-Y_1^{(j)}|=|Y_2^{(i)}-Y_2^{(j)}|~~\forall~i, j.\label{eq: E5}\tag{E5}
	\end{equation}
	We now claim that for $i, i^{'}\in\{1, 2, \cdots , n\}$, $\sum_{l=1}^{n}\ks(i/n,l/n)=\sum_{l=1}^{n} \ks(i^{'}/n,l/n)$ if and only if either $i^{'}=i$ or $i^{'}=n+1-i$. The `if' part 
	is easy to see; if $i^{'}=n+1-i$, the equality is obtained by observing that $\ks (i^{'}/n,j/n)=\ks(i/n, (n+1-j)/n)~\forall~j$ and then making a change of variable ($j\mapsto n+1-j$) in the summation. The 'only if' part can now be completed by observing that whenever $i<n+1-i$, $\sum_{l=1}^{n}\ks((i+1)/n,l/n)-\sum_{l=1}^{n}\ks(i/n,l/n)=e^{-\frac{i^2}{2n^2\sigma^2}}-e^{-\frac{(n-i)^2}{2n^2\sigma^2}}>0$, implying that
	$\sum_{l=1}^{n} \ks(i/n,l/n)$ is strictly increasing in $i$ whenever $i<n+1-i$.
	
	Using this, \eqref{eq: E4} implies that, for each $i$, we have either $Y^{(i)}_2=Y^{(i)}_1$ or $Y^{(i)}_2=1+\frac{1}{n}-Y^{(i)}_1$. Next, let $i$ be such that $Y_1^{(i)}=1/n$. We know that either $Y_2^{(i)}=Y_1^{(i)}$ or $Y_2^{(i)}=1+1/n-Y_1^{(i)}$. Suppose first that $Y_2^{(i)}=Y_1^{(i)}$. Now, take any $j\not=i$. We know $Y_2^{(j)}$ equals either $Y_1^{(j)}$ or $1+1/n-Y_1^{(j)}$. But then \eqref{eq: E5} rules out the possibility that $Y_2^{(j)}=1+1/n-Y_1^{(j)}$. Thus we have $Y_2^{(j)}=Y_1^{(j)}$ for all $j$. Similarly, if $Y_2^{(i)}=1+1/n-Y_1^{(i)}$, one can show that $Y_2^{(j)}=1+1/n-Y_1^{(j)}$ for all $j$. Thus we conclude that either $Y_2^{(j)}=Y_1^{(j)}~\forall~j$ or $Y_2^{(j)}=1+1/n-Y_1^{(j)}~\forall~j$. But this means that one coordinate of the observation vectors is either an increasing or a decreasing function of the other coordinate. \qed
\end{proof}

The following well-known result, which can be found in \cite{tsukahara2005semiparametric}, is crucial in our derivation of the limiting distributions of ${\widehat I}_{\sigma,n}(\Xvec)$ in both under the null and the alternative hypotheses.

\begin{customthm}{T1}[Weak convergence of copula process]
	\label{thm: cop conv}
	Let $\mathbf{X^{(1)}},\cdots,\mathbf{X^{(n)}}$ be independent observations on the random vector $\mathbf{X}$ with copula distribution $\cx$ and let $\cxn$ be the empirical copula based on $\mathbf{X^{(1)}},\cdots,\mathbf{X^{(n)}}$. If, for all $i=1,2,\cdots,d$, the $i^{\text{th}}$ partial derivatives $D_i\cx(u)$ of $\cx$ exist and are continuous, then the process $\sqrt{n}(\cxn-\cx)$ converges weakly in $l^{\infty}([0,1]^d)$ to the process $\mathbb{G}_{\cx}$ given by \vspace{-2ex}
	$$\mathbb{G}_{\cx}(u)=\mathbb{B}_{\cx}(u)-\sum_{i=1}^{d}D_i\cx(u)\mathbb{B}_{\cx}(u^{(i)}),\vspace{-2ex}$$
	where $\mathbb{B}_{\cx}$ is a $d$-dimensional Brownian bridge on $[0,1]^d$ with covariance function $\E[\mathbb{B}_{\cx}(u)\mathbb{B}_{\cx}(v)]=\cx(u)\wedge\cx(v)-\cx(u)\cx(v)$, and for each $i$, $u^{(i)}$ denotes the vector obtained from $u$ by replacing its all coordinates, except the $i^{\text{th}}$ th one, by $1$.
\end{customthm}

Let the distribution function of $\Xvec$ be denoted by $\mathbf{F}$ and its marginals by $F_1,F_2,\cdots,F_d$. With $\Xvec^{(i)}\!=\!(X_1^{(i)},X_2^{(i)}, \cdots X_d^{(i)}),\; 1\!\leq\! i\!\leq\! n$, denoting i.i.d.~observations from $\Xvec$, define vectors $\Zvec^{(i)}\!=\!(Z_1^{(i)},Z_2^{(i)},\allowbreak\cdots,Z_d^{(i)}),\; 1\!\leq\! i\!\leq\! n$, where $Z_j^{(i)}=F_j(X_j^{(i)})$. We will denote the empirical distribution based on $\Zvec^{(1)},\Zvec^{(2)},\cdots,\Zvec^{(n)}$ by $\mathsf{C}_{\zvec,n}$, and the empirical distribution function based on the $\Xvec^{(i)}$'s by $\hat{\mathbf{F}}$.

\begin{customlemma}{L2}
	\label{lem: 1}
	Assume that $\{P_n\}_{n\geq 1}$ is a sequence of distributions over $[0,1]^d$. Then
	\begin{enumerate}
		\item $\left|\gamma^2_{\ks}(\Pi_n,P_n)-\gamma^2_{\ks}(\Pi,P_n)\right|=\mathcal{O}(n^{-2})$
		\item $\left|\gamma^2_{\ks}(\mathsf{M}_n,P_n)-\gamma^2_{\ks}(\mathsf{M},P_n)\right|=\mathcal{O}(n^{-2}).$
	\end{enumerate}
\end{customlemma}

\begin{proof}
	We prove the first part only. The proof of the second part is similar.
	\begin{align*}
	\left|\gamma^2_{\ks}(\Pi_n,P_n)\!-\!\gamma^2_{\ks}(\Pi,P_n)\right|&\leq \left|\E_{_{(S,S^{'})\sim \Pi_n\otimes\Pi_n}}\!	[k_{\sigma}(S,S^{'})]-\E_{_{(S,S^{'})\sim \Pi\otimes\Pi}}\![k_{\sigma}(S,S^{'})]\right|\\
	&\quad+2\E_{_{T\sim P_n}}\left|\E_{_{S\sim \Pi_n}}\![k_{\sigma}(S,T)]-\E_{_{S\sim\Pi}}\![k_{\sigma}(S,T)]\right|.
	\end{align*}
	The first term on the right hand side of the above inequality is easily seen to be bounded above by \vspace{-1ex}
	$$\frac{1}{n^{2d}}\!\!\sum_{\substack{\mu=(i_1/n,\cdots,i_d/n)\\1\leq i_1,\cdots,i_d\leq n}}\!\,\sum_{\substack{\nu=(j_1/n,\cdots,j_d/n)\\1\leq j_1,\cdots,j_d\leq n}}\;\int\limits_{[\mu-1/n,\mu]}
	\int\limits_{[\nu-1/n,\nu]}\!\!\!\!\!\left| \ks(\mu,\nu)\!-\!\ks(\zeta,\eta)\right|d\zeta\,d\eta,\vspace{-1ex}$$
	where for any $u=(u_1, u_2,\cdots ,u_d)\in [0,1]^d$ and $\delta>0$, $[u-\delta, u]$ denotes the rectangle $[u_1-\delta, u_1]\times [u_2-\delta, u_2]\times \cdots \times [u_d-\delta, u_d]$. The last expression is clearly bounded above by \vspace{-1ex}
	$$\max_{\substack{\mu=(i_1/n,\cdots,i_d/n)\\1\leq i_1,\cdots,i_d\leq n}}\;\max_{\substack{\nu=(j_1/n,\cdots,j_d/n)\\1\leq j_1,\cdots,j_d\leq n}}\;\sup_{\zeta\in[\mu-1/n,\mu]}\;\sup_{\eta\in[\nu-1/n,\nu]}\left| \ks(\mu,\nu)-\ks(\zeta,\eta)\right|.\vspace{-1ex}$$
	Using Lemma 6 of \cite{poczos2012}, one can further deduce that the last expression is bounded above by $dn^{-2}$. \\
	Similar technique can be used for the second term to get the upper  bound \vspace{-1ex}
	$$2\,\E_{T\sim P_n}\; \max_{\substack{\mu=(i_1/n,i_2/n,\cdots,i_d/n)\\1\leq i_1,i_2,\cdots,i_d\leq n}}\;\sup_{\eta\in[\mu-1/n,\mu]}\left| \ks(\mu,T)-\ks(\eta,T)\right|\leq 2dn^{-2}.\vspace{-1ex}$$
	Combining these two bounds, we get  $\left|\gamma^2_{\ks}(\Pi_n,P_n)-\gamma^2_{\ks}(\Pi,P_n)\right|=\mathcal{O}(n^{-2})$. \qed
\end{proof}

\begin{customlemma}{L3}
	\label{lem: 2}
	$\gamma_{\ks}(\cxn,\mathsf{C}_{\zvec,n})\rightarrow 0$ almost surely as $n\rightarrow \infty$.
\end{customlemma}

\begin{proof}[{\bf Sketch of the proof}]
	Since the essential idea of the proof is contained in  \cite{poczos2012} [Appendix E], we only describe the two main steps.
	
	First, we use the definition of $\mathsf{C}_{\zvec,n}$ and Lemma 6 of \cite{poczos2012} to get the inequality
	$\gamma^2_{\ks}(\cxn,\mathsf{C}_{\zvec,n})\leq 2\sqrt{d}L\,\max\limits_{1\leq j\leq d}\;
	\sup\limits_{x\in\mathbb{R}}|\hat{F}_j(x)-F_j(x)|$, where $\hat{F}_1,\ldots,\hat{F}_d$ are the marginals of the empirical distribution $\hat{\mathbf{F}}$ based on
	$\Xvec^{(1)},\ldots,\Xvec^{(n)}$.
	
	Then using the above inequality and the Kiefer-Dvoretzky-Wolfowitz Theorem (see \cite{massart1990}, Page 1269), for any $\epsilon>0$, we get  $\pr\left[\gamma^2_{\ks}(\cxn,\mathsf{C}_{\zvec,n})>\epsilon\right]\leq 2d\exp\left(-\frac{n\epsilon^2}{2dL^2}\right)$. The result now follows from the Borel-Cantelli Lemma. \qed
\end{proof}

\vspace{0.025in}
The next lemma and its proof are based on the ideas in \cite{gretton2012} [Appendix A2].

\vspace{0.05in}
\begin{customlemma}{L4}
	\label{lem: 3}
	$\gamma_{\ks}(\mathsf{C}_{\zvec,n},\cx)\rightarrow 0$ almost surely as $n\rightarrow \infty$.
\end{customlemma}

\begin{proof}
	It is enough to prove $\E\left[\gamma_{\ks}(\mathsf{C}_{\zvec,n},\cx)\right]\leq \frac{2}{\sqrt{n}}\;
	\mbox { and }\;$\\
	$\pr\left[\gamma_{\ks}(\mathsf{C}_{\zvec,n},\cx)-\E[\gamma_{\ks}(\mathsf{C}_{\zvec,n},\cx)]>\epsilon\right]\leq \exp\left(\!-\frac{n\epsilon^2}{2}\right).$
	
	Denoting $\mathcal{F}$ to be the unit ball in the RKHS associated to the kernel $\ks$ on $\mathbb{R}^d$, one gets $\gamma_{\ks}(\mathsf{C}_{\zvec,n},\cx)= \sup_{f \in\mathcal{F}}\left|\frac{1}{n}\sum_{i=1}^{n}f(\Zvec^{(i)})-\E_{\Zvec\sim \cx}f(\Zvec)\right|$ (see \cite{sriperumbudur2010hilbert}).
	
	Letting $\Zvec^{*(i)}, \; 1\leq i\leq n$ to be i.i.d. with the same distribution as and independent of $\Zvec^{(i)},\; 1\leq i\leq n$ and $\delta_i,\; 1\leq i\leq n$ to be i.i.d. random variables taking values $\pm 1$ with equal probabilities, independent of the  $\Zvec^{(i)}, \Zvec^{*(i)}, \; 1\leq i\leq n$, it is easy to see that 
	\vspace{-10pt}
	\begin{align*}
	\E[\gamma_{\ks}(\mathsf{C}_{\zvec,n},\cx)]&=\E\left[\sup_{f \in \mathcal{F}}\left|\frac{1}{n}\sum_{i=1}^{n}f(\Zvec^{(i)})-\E_{\Zvec\sim \cx}f(\Zvec)\right|\right]\\
	&\leq \E\left[\sup_{f \in \mathcal{F}}\left|\frac{1}{n}\sum_{i=1}^{n}f(\Zvec^{(i)})-\frac{1}{n}\sum_{i=1}^{n}f(\Zvec^{*(i)})\right|\right]\\
	&=\E\left[\sup_{f \in \mathcal{F}}\left|\frac{1}{n}\sum_{i=1}^{n}\delta_i\left(f(\Zvec^{(i)})-f(\Zvec^{*(i)})\right)\right|\right]	\underset{(a)}{\leq} \frac{2}{\sqrt{n}}.
	\vspace{-10pt}
	\end{align*}
	For the last inequality (a), we used a well-known result referred to as ``Bound on Rademacher Complexity" (see \cite{bartlett2003}, Page 478) .
	
	We next calculate the upper bound of change in magnitude due to change in a particular coordinate. Consider $\gamma_{\ks}(\mathsf{C}_{\zvec,n},\cx)$ as a function of $\Zvec^{(i)}$. It is easy to verify that changing any coordinate $\Zvec^{(i)}$, the change in $\gamma_{\ks}(\mathsf{C}_{\zvec,n},\cx)$ will be at most $2n^{-1}$. We use now the well-known McDiarmid's inequality (see \cite{mcdiarmid1989}, Page 149) to get
	\vspace{-10pt}
	$$\pr\left[\gamma_{\ks}(\mathsf{C}_{\zvec,n},\cx)-\E[\gamma_{\ks}(\mathsf{C}_{\zvec,n},\cx)]>\epsilon\right]\leq \exp\left(-\frac{2\epsilon^2}{n.(2/n)^2}\right)=\exp\left(-\frac{n\epsilon^2}{2}\right). \qed $$
\end{proof}
\vspace{-10pt}

\begin{customlemma}{L5}
	\label{lem: asymp conv}
	Suppose that the assumptions of Theorem \ref{thm: cop conv} hold. Then, we have the following results.
	\begin{itemize}
		\item[]If $\cx\neq\Pi$,
	$\sqrt{n}(\gamma_{\ks}^2(\cxn,\Pi)-\gamma_{\ks}^2(\cx,\Pi))\overset{\mathcal{L}}{\rightarrow}\mathcal{N}(0,\delta_0^2),$
	where\\ $\delta^2_0=$ $4\displaystyle\int_{\mathbb{R}^d}\displaystyle\int_{\mathbb{R}^d}{\hspace{-0.1in}g(u)}g(v)\,\E[\,d\mathbb{G}_{\cx}(u)\,d\mathbb{G}_{\cx}(v)]$ and $g(u)=\displaystyle\int_{[0,1]^d} {\hspace{-0.25in} \ks(u,v)\,d(\cx-\Pi)(v)}$.
	\item[]If $\cx=\Pi$,
	$n\gamma_{\ks}^2(\cxn,\Pi)\overset{\mathcal{L}}{\rightarrow}\displaystyle\int_{\mathbb{R}^d}\displaystyle\int_{\mathbb{R}^d}\ks(u,v)\,d\mathbb{G}_\Pi(u)\,d\mathbb{G}_\Pi(v).$
	\end{itemize}
\end{customlemma}

\begin{proof} \underline{When $\cx\neq\Pi$:~}
	Denoting $\mathcal{D}([0,1]^d)$ to be the space of right continuous real valued uniformly bounded functions on $[0,1]^d$ with left limits, equipped with max-sup norm, one can easily verify that the function $\psi(\mathsf{D})=\gamma^2_{\ks}(\mathsf{D},\Pi)$ on $\mathcal{D}([0,1]^d)$ is Hadamard-differentiable and the derivative at $\cx$ is given by $$\psi^{'}_{\cx}(\mathsf{D})=2\int_{[0,1]^d}\int_{[0,1]^d}\ks(u,w)\,d(\cx-\Pi)(w)\,d\mathsf{D}(u).$$
	To prove this, consider a real sequence $\{t_n\}$  converging to 0 and a $\mathcal{D}([0,1]^d)$-valued sequence $\{D_n\}$ converging to $D\in\mathcal{D}([0,1]^d)$ such that $\cx+t_n D_n\in\mathcal{D}([0,1]^d)$. For any $D\in\mathcal{D}([0,1]^d)$, define $\mu_D(w)=\intc\kss(u,w)\,dD(u)~\forall w\in\Rd$. 
	Then
	\begin{align}
	&\frac{\varphi(\cx+t_nD_n)-\varphi(\cx)}{t_n}\nonumber\\
	&=\left(\frac{1}{\sigma}\sqrt{\frac{2}{\pi}}\right)^d\frac{1}{t_n}\int_{\Rd}\left(\mu_{\cx}(w)+t_n\mu_{D_n}(w)-\mu_\Pi(w)\right)^2\,dw\nonumber\\
	&\qquad-\left(\frac{1}{\sigma}\sqrt{\frac{2}{\pi}}\right)^d\frac{1}{t_n}\int_{\Rd}\left(\mu_{\cx}(w)-\mu_\Pi(w)\right)^2\,dw\nonumber\\
	&=\left(\frac{1}{\sigma}\sqrt{\frac{2}{\pi}}\right)^d\frac{1}{t_n}\int_{\Rd}t_n\mu_{D_n}(w)\left(2\mu_{\cx}(w)+t_n\mu_{D_n}(w)-2\mu_\Pi(w)\right)\,dw\nonumber\\
	&=\left(\frac{1}{\sigma}\sqrt{\frac{2}{\pi}}\right)^d\left[2\int_{\Rd}\mu_{D_n}(w)\left(\mu_{\cx}(w)-\mu_\Pi(w)\right)\,dw+t_n\int_{\Rd}\mu^2_{D_n}(w)\,dw\right]\label{eq: hadamard}
	\end{align}
	Now, using the fact $\int_{\Rd}\left(\frac{1}{\sigma}\sqrt{\frac{2}{\pi}}\right)^d\kss(u,w)\kss(v,w)\,dw=\ks(u,v)$, it is quite straight-forward to check that $\left(\frac{1}{\sigma}\sqrt{\frac{2}{\pi}}\right)^d\int_{\Rd}\mu_{D_n}(w)\left(\mu_{\cx}(w)-\mu_\Pi(w)\right)\,dw=$\\
	$\intc\intc\ks(u,v)\,dD_n(u)\,d(\cx-\Pi)(v)$. From this identity and Equation \eqref{eq: hadamard}, we get
	\begin{align*}
	\psi^{'}_{\cx}(\mathsf{D})&=\lim_{n\to\infty}
	\frac{\varphi(\cx+t_nD_n)-\varphi(\cx)}{t_n}\\&=2\intc\intc\ks(u,v)\,dD(u)\,d(\cx-\Pi)(v).
	\end{align*}
	
	This Lemma then follows easily from Theorem \ref{thm: cop conv} and the functional delta method. The only thing that one needs to verify is that $\psi^{'}_{\cx}(\mathbb{G}_{\cx})$ is a normally distributed with 0 mean and variance $\delta^2_0=4\int_{[0,1]^d}\int_{[0,1]^d}g(u)g(v)\,\E[\,d\mathbb{G}_{\cx(}u)\,d\mathbb{G}_{\cx}(v)]$ where $g(u)=\int_{[0,1]^d}\ks(u,w)\,d(\cx-\Pi)(w)$. But this is straightforward from the formula for the derivative $\psi^{'}_{\cx}$.
	
	\noindent\underline{When $\cx=\Pi$:~}
	Clearly the map $\mathsf{D}\rightarrow\left(\frac{1}{\sigma}\sqrt{\frac{2}{\pi}}\right)^d\!\int_{\mathbb{R}^d}\left(\int_{[0,1]^d}k_{\frac{\sigma}{\sqrt{2}}}(u,v)\,d\mathsf{D}(u)\right)^2\!dv$ from $\mathcal{D}([0,1]^d)$ to $\mathbb{R}$ is  continuous. So, the fact that $\sqrt{n}(\cxn-\Pi)\rightarrow\mathbb{G}_{\Pi}$ and the continuous mapping theorem gives\vspace{-1ex}
	\begin{align*}
	n\gamma^2_{\ks}(\cxn,\Pi)&\overset{\mathcal{L}}{\longrightarrow}\left(\frac{1}{\sigma}\sqrt{\frac{2}{\pi}}\right)^d.\int_{\mathbb{R}^d}\left(\int_{[0,1]^d}k_{\frac{\sigma}{\sqrt{2}}}(u,v)\,d\mathbb{G}_\Pi(u)\right)^2\,dv\\
	&=\int_{[0,1]^d}\int_{[0,1]^d}\ks(u,v)\,d\mathbb{G}_\Pi(u)\,d\mathbb{G}_\Pi(v). \quad\quad\quad\quad\quad\quad\quad\quad~~~~~~~~~\qed
	\end{align*}
\end{proof}

\begin{proof}[{\bf Proof of Theorem 5}]
	\hspace{0pt}\\
	\underline{When $\cx\neq\Pi$:~}
	Write $\sqrt{n}(\secb\!-\!I^2_\sigma)\!=\!\sqrt{n}\left(\frac{\gamma^2_{\ks}(\cxn,\Pi_n)}
	{\gamma^2_{\ks}(\mathsf{M}_n,\Pi_n)}\!-\!\frac{\gamma^2_{\ks}(\cx,\Pi)}
	{\gamma^2_{\ks}(\mathsf{M},\Pi)}\right)$ as\\
	$A_{1,n}+A_{2,n}+A_{3,n}$, where \vspace{-10pt}
	\begin{align*}
	&A_{1,n}=\sqrt{n}\left(\frac{\gamma^2_{\ks}(\cxn,\Pi_n)}{\gamma^2_{\ks}(\mathsf{M}_n,\Pi_n)}-\frac{\gamma^2_{\ks}(\cxn,\Pi_n)}{\gamma^2_{\ks}(\mathsf{M},\Pi)}\right),\\ &A_{2,n}=\sqrt{n}\left(\frac{\gamma^2_{\ks}(\cxn,\Pi_n)}{\gamma^2_{\ks}(\mathsf{M},\Pi)}-\frac{\gamma^2_{\ks}(\cxn,\Pi)}{\gamma^2_{\ks}(\mathsf{M},\Pi)}\right)\!\\
	&\mbox { and } A_{3,n}=\sqrt{n}\left(\frac{\gamma^2_{\ks}(\cxn,\Pi)}{\gamma^2_{\ks}(\mathsf{M},\Pi)}-\frac{\gamma^2_{\ks}(\cx,\Pi)}{\gamma^2_{\ks}(\mathsf{M},\Pi)}\right).\vspace{-10pt}
	\end{align*}
	Clearly $A_{2,n} \rightarrow 0$ almost surely by Lemma \ref{lem: 1}. The same is true of $A_{1,n}$ as well, once again by Lemma \ref{lem: 1}, because it is bounded above by 	\vspace{-10pt}
	$$\underbrace{\frac{\gamma^2_{\ks}\!(\cxn,\!\Pi_n)}{\gamma^2_{\ks}\!(\mathsf{M}_n,\!\Pi_n)
			\gamma^2_{\ks}\!(\mathsf{M},\!\Pi)}}_{\text{bounded sequence}}\left(\underbrace{\sqrt{n}\left|\gamma^2_{\ks}\!(\mathsf{M}_n,\!\Pi_n)\!-\!\gamma^2_{\ks}\!(\mathsf{M},\!\Pi_n)\right|}_{\text{goes to 0}}+\underbrace{\sqrt{n}\left|\gamma^2_{\ks}\!(\mathsf{M},\!\Pi_n)\!-\!\gamma^2_{\ks}\!(\mathsf{M},\!\Pi)\right|}_{\text{goes to 0}}\right).
	\vspace{-10pt} $$
	Therefore, using Lemma \ref{lem: asymp conv}, we can conclude that $\sqrt{n}(\secb-\scb)\overset{\mathcal{L}}{\rightarrow}\mathcal{N}(0,\csd^{-2}.\delta^2_0)$.  Now, applying the delta method, one gets 
	$\sqrt{n}(\ecb-\cb)\overset{\mathcal{L}}{\rightarrow}\mathcal{N}(0,\delta^2)$, where $\delta^2
	=\csd^{-2}I^{-2}_\sigma(\Xvec)\int_{[0,1]^d}\int_{[0,1]^d}\!g(u)g(v)\,\E[\,d\mathbb{G}_{\cx}(u)\,d\mathbb{G}_{\cx}(v)].$\\
	
	\noindent\underline{When $\cx=\Pi$:~}
	As a consequence of the Lemma L5 and Lemma \ref{lem: 1}, under null hypothesis and assumptions of Theorem \ref{thm: cop conv}, we have \vspace{-5pt}
	\begin{align*}
	n\gamma^2_{\ks}(\cxn,\Pi_n)&=n\gamma^2_{\ks}(\cxn,\Pi)+\left(n\gamma^2_{\ks}(\cxn,\Pi_n)-n\gamma^2_{\ks}(\cxn,\Pi)\right)\\ &\overset{\mathcal{L}}{\longrightarrow}\;\int_{[0,1]^d}\int_{[0,1]^d}\ks(u,v)\,d\mathbb{G}_\Pi(u)\,d\mathbb{G}_\Pi(v).\vspace{-10pt}
	\end{align*}
	It is enough to show that $\int_{[0,1]^d}\int_{[0,1]^d}\ks(u,v)\,d\mathbb{G}_\Pi(u)\,d\mathbb{G}_\Pi(v)\overset{\mathcal{L}}{=}\sum_{i=1}^{\infty}\alpha_iZ_i^2,$ for some $\alpha_i>0$ and $Z_i\overset{\text{i.i.d.}}{\sim}\mathcal{N}(0,1)$. The actual result will follow putting $\lambda_i=\alpha_i\csd$.
	
	To this end, we define $X(w):=\int_{[0,1]^d}k_{\frac{\sigma}{\sqrt{2}}}(u,w)\,d\mathbb{G}_\Pi(u),~\forall w\in\mathbb{R}^d$. So, $\{X(w):w\in\mathbb{R}^d\}$ is then a zero-mean continuous path Gaussian process. This implies that
	$$\int_{\mathbb{R}^d}\left(\int_{[0,1]^d}k_{\frac{\sigma}{\sqrt{2}}}(u,w)\,d\mathbb{G}_\Pi(u)\right)^2\!dw\;=\int_{\mathbb{R}^d}(X(w))^2\,dw\;\overset{\mathcal{L}}{=}\;\sum_{i=1}^{\infty}\beta_iZ_i^2.\vspace{-1.5ex} $$
	Now, using the fact that $\left(\frac{1}{\sigma}\sqrt{\frac{2}{\pi}}\right)^d\int_{\mathbb{R}^d}k_{\frac{\sigma}{\sqrt{2}}}(u,w)k_{\frac{\sigma}{\sqrt{2}}}(v,w)dw=\ks(u,v)$, one can easily see that the last equality yields the desired result\vspace{-1.5ex}
	$$\int_{[0,1]^d}\int_{[0,1]^d}\ks(u,v)\,d\mathbb{G}_\Pi(u)\,d\mathbb{G}_\Pi(v)\;\overset{\mathcal{L}}{=}\;\sum_{i=1}^{\infty}\alpha_iZ_i^2,\;\; \mbox { with }  \alpha_i=\left(\frac{1}{\sigma}\sqrt{\frac{2}{\pi}}\right)^d\!\!\beta_i. ~~~\qed \vspace{-1.5ex}$$
\end{proof}

\begin{proof}[{\bf Proof of Theorem 6}] Triangle inequality and $|a-b|^2\leq|a^2-b^2|$ for $a,b\geq 0$ gives \vspace{-5pt}
	$$\left|\gamma_{\ks}\!(\mathsf{M}_n,\!\Pi_n)\!-\!\gamma_{\ks}\!(\mathsf{M},\!\Pi)\right|\!\leq\!
	\left|\gamma^2_{\ks}\!(\mathsf{M}_n,\!\Pi_n)\!-\!\gamma^2_{\ks}\!(\mathsf{M}_n,\!\Pi)\right|^{\!\frac{1}{2}}+
	\left|\gamma^2_{\ks}\!(\mathsf{M}_n,\!\Pi)\!-\!\gamma^2_{\ks}\!(\mathsf{M},\!\Pi)\right|^{\!\frac{1}{2}}. $$
	Using Lemma \ref{lem: 1}, we get $\lim_{n\rightarrow\infty}\left|\gamma_{\ks}(\mathsf{M}_n,\Pi_n)-\gamma_{\ks}(\mathsf{M},\Pi)\right|=0$ a.s.. Using again the same inequalities and the fact that $\gamma_{\ks}$ is a metric, one gets
	\vspace{-5pt}
	\begin{align*}
	\left|\gamma_{\ks}\!(\cxn,\!\Pi_n)\!-\!\gamma_{\ks}\!(\cx,\!\Pi)\right|
	\!&\leq\!\left|\gamma^2_{\ks}\!(\cxn,\!\Pi_n)\!-\!\gamma^2_{\ks}\!(\cxn,\!\Pi)\right|^{\!\frac{1}{2}}\\
	&\quad+\gamma_{\ks}\!(\cxn,\!\mathsf{C}_{\zvec,n})+\gamma_{\ks}\!(\mathsf{C}_{\zvec,n},\!\cx).
	\end{align*}
	Again, using Lemmas \ref{lem: 1}, \ref{lem: 2} and \ref{lem: 3}, we get $\left|\gamma_{\ks}(\cxn,\Pi_n)-\gamma_{\ks}(\cx,\Pi)\right|\stackrel{a.s.}{\rightarrow} 0$  as $n\rightarrow 0$, and as a consequence, we conclude that as $n\rightarrow\infty$,
	\vspace{-5pt}
	$$\ecb=\frac{\gamma_{\ks}(\cxn,\Pi_n)}{\gamma_{\ks}(\mathsf{M}_n,\Pi_n)}\;\rightarrow\;\frac{\gamma_{\ks}(\cx,\Pi)}{\gamma_{\ks}(\mathsf{M},\Pi)}=\cb \mbox { almost surely} \quad\quad\quad\quad~~~~~~~~~~~\qed.$$ 
\end{proof}

\begin{customlemma}{L6}
	\label{lem: sig n}
	Let $P_n$ and $Q_n$ be sequence of  probability distribution over $[0,1]^d$. Let $\sigma_n$ be a sequence of positive real numbers that converges to $\sigma_0>0$. Then as $n\to\infty$, $|\gamma_{k_{\sigma_n}}^2(P_n,Q_n)-\gamma_{k_{\sigma_0}}^2(P_n,Q_n)|\to 0.$
\end{customlemma}

\begin{proof}
	First we observe that
	\vspace{-5pt}
	\begin{align*}
	|\gamma_{k_{\sigma_n}}^2(P_n,Q_n)-\gamma_{k_{\sigma_0}}^2(P_n,Q_n)|&\leq \E_{(\Svec,\Svec^{'})\sim P_n\otimes P_n}|k_{\sigma_n}(\Svec,\Svec^{'})-k_{\sigma_0}(\Svec,\Svec^{'})|\\
	&\quad+2\E_{(\Svec,\Tvec)\sim P_n\otimes Q_n}|k_{\sigma_n}(\Svec,\Tvec)-k_{\sigma_0}(\Svec,\Tvec)|\\
	&\quad+\E_{(\Tvec,\Tvec^{'})\sim Q_n\otimes Q_n}|k_{\sigma_n}(\Svec,\Tvec^{'})-k_{\sigma_0}(\Tvec,\Tvec^{'})|.
	\end{align*}
	Applying Lemma 6 of \cite{poczos2012}, one can get an upper bound of $|k_{\sigma_n}(\Svec,\Tvec)$ $-k_{\sigma_0}(\Svec,\Tvec)|$ in the following way
	\vspace{-5pt}
	\begin{align*}
	|k_{\sigma_n}(\Svec,\Tvec)-k_{\sigma_0}(\Svec,\Tvec)|&\leq L\left\|\frac{\Svec}{\sigma_n}-\frac{\Svec}{\sigma_0}\right\|+L\left\|\frac{\Tvec}{\sigma_n}-\frac{\Tvec}{\sigma_0}\right\|
	\leq 2L\sqrt{d}\left|\frac{1}{\sigma_n}-\frac{1}{\sigma_0}\right|,
	\end{align*}
	where $L$ is a constant. Thus we can conclude that
	$|\gamma_{k_{\sigma_n}}^2(P_n,Q_n)-\gamma_{k_{\sigma_0}}^2(P_n,Q_n)|\leq 8L\sqrt{d}\left|\frac{1}{\sigma_n}-\frac{1}{\sigma_0}\right|.$
	This completes the proof. \qed
\end{proof}

\begin{proof}[{\bf Proof of Theorem 7}] Note that
	\vspace{-3pt}
	\begin{align*}
	|\gamma^2_{\sigma_n}(\cxn,\Pi_n)-\gamma^2_{\sigma_0}(\cx,\Pi)|&\leq |\gamma^2_{\sigma_n}(\cxn,\Pi_n)-\gamma^2_{\sigma_0}(\cxn,\Pi_n)|\\
	&\quad+|\gamma^2_{\sigma_0}(\cxn,\Pi_n)-\gamma^2_{\sigma_0}(\cx,\Pi)|.
	\end{align*}
	As $n\to\infty$, the first term in the right hand side goes to 0 due to Lemma \ref{lem: sig n} and the second term almost surely converges to 0 due to Theorem 6. Similarly, one can show that $|\gamma^2_{\sigma_n}(\mathsf{M}_n,\Pi_n)-\gamma^2_{\sigma_0}(\mathsf{M},\Pi)|\to 0$ as $n\to\infty$. This implies $\widehat{I}_{\sigma_n,n}(\Xvec)\to I_{\sigma_0}(\Xvec)$ almost surely. Because of the fact that $I_{\sigma_0}(\Xvec)=0$ if and only if the coordinates of $\Xvec$ are independent, test of independence based on the statistic $\widehat{I}_{\sigma_n,n}(\Xvec)$ is consistent. \qed
\end{proof}

\begin{proof}[{\bf Proof of Theorem 8}]
	From Theorem (7), it follows that $\widehat{I}_{\sigma^{(i)},n}(\Xvec)$'s are consistent test statistics for all $i=1,2,\cdots,m$. Since $m$ is finite, by the virtue of the definition of $T_{max}$ and $T_{sum}$, they converge to 0 almost surely if and only if the coordinates of $\Xvec$ are independent. Otherwise, they converge to positive quantities. This property makes the resulting tests consistent. Again, under the alternative hypothesis, for any $i$, the p-value $p_i$ corresponding to the test statistic $\widehat{I}_{\sigma^{(i)},n}(\Xvec)$ converges to zero almost surely. So, for sufficiently large $n$, almost surely, there would exist at least one $i$ such that $p_i$ is less than $\alpha/m$, which makes the set $\{i:p_{(i)}<\alpha/m\}$ non-empty. Thus the power of the test based on FDR tends to be $1$ as sample size tends to infinity. \qed
\end{proof}

\bibliographystyle{spbasic}      
\bibliography{References}

\end{document}